\documentclass[a4paper,12pt]{article}

\usepackage{vmargin}            
\usepackage{fancyheadings}      
\usepackage{calc}               
\usepackage{fancybox}           
\usepackage{theorem}

\usepackage{amsmath}    
\usepackage{amssymb}    
\usepackage{amsfonts}   

\setmarginsrb{2cm}{2.15cm}{2cm}{2.15cm}{1mm}{6mm}{0mm}{10mm}

\newcommand{\f}[2]{\ensuremath{\mathchoice%
        {\dfrac{#1}{#2}}
       {\dfrac{#1}{#2}}
        {\frac{#1}{#2}}
        {\frac{#1}{#2}}
        }}

\newcommand{\dd}{\ensuremath{\text{d}}}

\def\d{{\rm d}}

\newcommand{\Prod}{\ensuremath{\prod\limits}}

\renewcommand{\emptyset}{\ensuremath{\varnothing}}

\newcommand{\N}{\ensuremath{\mathbb{N}}}
\newcommand{\Z}{\ensuremath{\mathbb{Z}}}

\newcommand{\R}{\ensuremath{\mathbb{R}}}
\newcommand{\T}{\ensuremath{\mathbb{T}}}
\newcommand{\C}{\ensuremath{\mathbb{C}}}

\renewcommand{\dim}{\ensuremath{\mathop{\rm dim\,}\nolimits}}

\newcommand{\Rez}{\ensuremath{\mathop{\text{Re}\,}\nolimits}}
\newcommand{\Imz}{\ensuremath{\mathop{\text{Im}\,}\nolimits}}

\renewcommand{\Re}{\ensuremath{\Rez}}
\renewcommand{\Im}{\ensuremath{\Imz}}

\renewcommand{\le}{\ensuremath{\leqslant}}
\renewcommand{\ge}{\ensuremath{\geqslant}}

\newcommand{\scal}[2]{\ensuremath{\langle #1,#2 \rangle}}
\theoremstyle{plain}
\newtheorem{thm}{Theorem}[section]
\newtheorem{lem}[thm]{Lemma}
\newtheorem{cor}[thm]{Corollary}
\newtheorem{prop}[thm]{Proposition}
\newtheorem{defi}[thm]{Definition}
\newtheorem{sch}[thm]{Scholie}
\newcommand{\proof}{\noindent {\it Proof.}\;}

\newcommand{\Eproof}{\hfill{{\large $\square$}}}

\newcommand\FQ{\mathcal{F}_{q}\left(H_{\C}\right)}

\newcommand\HC{H_{\C}}

\newcommand\id{\text{id}}

\begin{document}
\pagestyle{empty}
\renewcommand{\thefootnote}{}
~ \vspace{2cm}
\begin{center}
{ \Large \bf Asymptotic matricial models and QWEP property for $q$-Araki-Woods algebras}\\
\vspace{4 cm}
Alexandre Nou\\
Universit\'e de Franche-Comt\'e - Besancon\\
U.F.R des Sciences et Techniques\\
D\'epartement de Math\'ematiques\\
16 route de Gray - 25030 Besancon Cedex\\
nou@math.univ-fcomte.fr
\end{center}
\vspace{4 cm} {\bf Abstract}.{ Using Speicher central limit Theorem
we provide Hiai's q-Araki-Woods
  von Neumann algebras with nice asymptotic matricial models. Then, we use this model and an elaborated
  ultraproduct procedure,
   to show that all q-Araki-Woods von Neumann algebras are QWEP.}
\footnotetext[1]{AMS classification: 46L65, 46L53 \hfill Keywords:
QWEP, deformation, matricial models, ultraproduct}

\newpage
\pagestyle{plain}

\section{Introduction}
 Recall that a $C^{*}$-algebra has the weak expectation property (in short WEP) if
the canonical inclusion from $A$ into $A^{**}$ factorizes completely
contractively through some $B(H)$ ($H$ Hilbert). A $C^{*}$-algebra
is QWEP if it is a quotient by a closed ideal of an algebra with the
WEP. The notion of QWEP was introduced by Kirchberg in \cite{Kir}.
Since then, it became an important notion in the theory of
$C^{*}$-algebras. Very recently, Pisier and Shlyakhtenko \cite{PS}
proved that Shlyakhtenko's free quasi-free factors are QWEP. This
result plays an important role in their work on the operator space
Grothendieck Theorem, as well as in the subsequent related works
\cite{P1} and \cite{X}. On the other hand, in his paper \cite{J} on
the embedding of Pisier's operator Hilbertian space $OH$ and the
projection constant of $OH_{n}$, Junge used QWEP in a crucial way.

\medskip

Hiai \cite{Hi} introduced the so-called $q$-Araki-Woods algebras.
Let $-1<q<1$, and let $H_{\R}$ be a real Hilbert space and
$(U_{t})_{t \in \R}$ an orthogonal group on  $H_{\R}$. Let
$\Gamma_{q}(H_{\R}, (U_{t})_{t \in \R})$ denote the associated
$q$-Araki-Woods algebra. These algebras are generalizations of both
Shlyakhtenko's free quasi-free factors (for $q=0$), and Bo${\rm
\dot{z}}$ejko and Speicher's $q$-Gaussian algebras (for $(U_{t})_{t
\in \R}$ trivial). In this paper we prove that $\Gamma_{q}(H_{\R},
(U_{t})_{t \in \R})$ is QWEP. This is an extension of
Pisier-Shlyakthenko's result for the free quasi-free factor (with
$(U_{t})_{t \in \R}$ almost periodic), already quoted above.

\medskip

In the first two sections below we recall some general background on
$q$-Araki-Woods algebras and we give a proof of our main result in
the particular case of Bo${\rm \dot{z}}$ejko and Speicher's
$q$-Gaussian algebras $\Gamma_{q}(H_{\R})$. The proof relies on an
asymptotic random matrix model for standard $q$-Gaussians. The
existence of such a model goes back to Speicher's central limit
Theorem for mixed commuting/anti-commuting non-commutative random
variables (see \cite{Sp}). Alternatively, one can also use the
Gaussian random matrix model given by \'Sniady in \cite{Sn}. Notice
that the matrices arising from Speicher's central limit Theorem may
not be uniformly bounded in norm. Therefore, we have to cut them off
in order to define a homomorphism from a dense subalgebra of
$\Gamma_{q}(H_{\R})$ into an ultraproduct of matricial algebras. In
this tracial framework it can be shown quite easily that this
homomorphism extends to an isometric $*$-homomorphism of von Neumann
algebras, simply because it is trace preserving. Thus
$\Gamma_{q}(H_{\R})$ can be seen as a (necessarily completely
complemented) subalgebra of an ultraproduct of matricial algebras.
This solves the problem in the tracial case.

Moreover, in this (relatively) simple situation, we are able to
extend the result to the $C^*$-algebra generated by all
$q$-Gaussians, $C^{*}_{q}(H_{\R})$. Indeed, using the
ultracontractivity of the $q$-Ornstein Uhlenbeck semi-group (see
\cite{B}) we establish that $C^{*}_{q}(H_{\R})$ is "weakly ucp
complemented" in $\Gamma_{q}(H_{\R})$. This last fact, combined with
the QWEP of $\Gamma_{q}(H_{\R})$, implies that $C^{*}_{q}(H_{\R})$
is also QWEP.

\medskip

In the remaining of the paper we adapt the proof of section
\ref{section2} to the more general type-$I\!I\!I$ $q$-Araki-Woods
algebras. In section \ref{section3} we start by recalling Raynaud's
construction of the von Neumann algebra's ultraproduct when algebras
are equipped with non-tracial states (see \cite{Ray}). Then, we give
some general conditions in order to define an embedding into such an
ultraproduct, whose image is of a state preserving conditional
expectation.

In section \ref{section4} we define a twisted Baby Fock model, to
which we apply Speicher's central limit Theorem. This provides us
with an asymptotic random matrix model for (finite dimensional)
$q$-Araki Woods algebras, generalizing the asymptotic model already
introduced by Speicher and used by Biane in \cite{Bi}. Using this
asymptotic model, we then define an algebraic $*$-homomorphism from
a dense subalgebra of $\Gamma_{q}(H_{\R}, (U_{t})_{t \in \R})$ into
a von Neumann ultraproduct of finite dimensional $C^{*}$-algebras.
Notice that the cut off argument requires some extra work (compare
the proofs of Lemma \ref{pres} and Lemma \ref{pres2}), for instance
we
 need to use our knowledge of the modular theory at the Baby Fock level
 to conclude. We then
apply the general results of section \ref{section3} (Theorem
\ref{qwep}) to extend this algebraic $*$-homomorphism into a
$*$-isomorphism from $\Gamma_{q}(H_{\R}, (U_{t})_{t \in \R})$ to the
von Neumann algebra's ultraproduct, whose image is completely
complemented. This allows us to show that $\Gamma_{q}(H_{\R},
(U_{t})_{t \in \R})$ is QWEP for $H_{\R}$ finite dimensional (see
Theorem \ref{qwepf}). It implies, by inductive limit, that
$\Gamma_{q}(H_{\R}, (U_{t})_{t \in \R})$ is QWEP when $(U_{t})_{t
\in \R}$ is almost periodic (see Corollary \ref{qwepa}).

In the last section, we consider a general algebra
$\Gamma_{q}(H_{\R}, (U_{t})_{t \in \R})$. We use a discretization
procedure on  the unitary group $(U_{t})_{t \in \R}$ in order to
approach $\Gamma_{q}(H_{\R}, (U_{t})_{t \in \R})$ by almost periodic
q-Araki-Woods algebras. We then apply the general results of section
\ref{section3} and, we recover the general algebra as a complemented
subalgebra of the ultraproduct of the discretized ones (see Theorem
\ref{qwepc}). From this last fact follows the QWEP of
$\Gamma_{q}(H_{\R}, (U_{t})_{t \in \R})$. However we were unable to
establish the corresponding result for the $C^{*}$-algebra
$C_{q}^{*}(H_{\R}, (U_{t})_{t \in \R})$. Indeed, if $(U_{t})_{t \in
\R}$ is not trivial then the ultracontractivity of the
$q$-Ornstein-Uhlenbeck semi-group never holds in any
right-neighborhood of zero (see \cite{Hi}).

We highlight that the modular theory on the twisted Baby Fock
algebras, on their ultraproduct, and on the $q$-Araki Woods
algebras, are crucial tools in order to overcome the difficulties
arising in the non-tracial case.

After the completion of this work, Marius Junge informed us that he
had obtained our main result using his proof of the non-commutative
L$^{1}$-Khintchine inequalities for $q$-Araki-Woods algebras.
Junge's approach is slightly different but its main steps are the
same as ours: the proof uses in a crucial way Speicher's central
limit Theorem, an ultraproduct argument and modular theory.

\section{Preliminaries}\label{section1}

\subsection{$q$-Araki-Woods algebras}\label{subesc1}
We mainly follow the notations used in \cite{Sh}, \cite{Hi} and
\cite{Nou}. Let $H_{\R}$ be a real Hilbert space and $(U_{t})_{t \in
\R}$ be a strongly continuous group of orthogonal transformations on
$H_{\R}$. We denote by $H_{\C}$ the complexification of $H_{\R}$ and
still by $(U_{t})_{t \in \R}$ its extension to a group of unitaries
on $H_{\C}$. Let $A$ be the (unbounded) non degenerate positive
infinitesimal generator of $(U_{t})_{t \in \R}$.
$$U_{t}=A^{it}\qquad\text{for all\;} t \in \R$$
A new scalar
 product $\scal{\;.\;}{\;.\;}_{U}$ is defined on $H_{\C}$ by the
 following relation:
$$\scal{\xi}{\eta}_{U}=\scal{2A(1+A)^{-1}\xi}{\eta}$$
We denote by $H$ the completion of $H_{\C}$ with respect to this new
scalar product. For $q \in (-1,1)$ we consider the $q$-Fock space
associated with $H$ and given by:
$$\mathcal{F}_{q}(H)=\C\Omega \bigoplus_{n\ge 1}H^{\otimes
n}$$ where $H^{\otimes n}$ is equipped with Bo${\rm \dot{z}}$ejko
and Speicher's $q$-scalar product (see \cite{BS1}). The usual
creation and annihilation operators on $\mathcal{F}_{q}(H)$ are
denoted respectively by $a^{*}$ and $a$ (see \cite{BS1}). For $f\in
H_{\R}$, $G(f)$, the $q$-Gaussian operator associated to $f$, is by
definition:
$$G(f)=a^{*}(f)+a(f)\in B\left(\mathcal{F}_{q}(H)\right)$$
The von Neumann algebra that they generate in
$B\left(\mathcal{F}_{q}(H)\right)$ is the so-called $q$-Araki-Woods
algebra: $\Gamma_{q}(H_{\R},(U_{t})_{t \in \R})$. The
$q$-Araki-Woods algebra is equipped with a faithful normal state
$\varphi$ which is the expectation on the vacuum vector $\Omega$. We
denote by $W$ the Wick product ; it is the inverse of the mapping:
\begin{equation*}
\begin{split}
  \Gamma_{q}(H_{\R},(U_{t})_{t \in \R}) & \longrightarrow
\Gamma_{q}(H_{\R},(U_{t})_{t \in \R})\Omega\\
X & \mapsto X\Omega
\end{split}
\end{equation*}
Recall that $\Gamma_{q}(H_{\R},(U_{t})_{t \in \R}) \subset
B\left(\mathcal{F}_{q}\left(H\right)\right)$ is the GNS
representation of $(\Gamma, \varphi)$. The modular theory relative
to the state $\varphi$ was computed in the papers \cite{Hi} and
\cite{Sh}. We now briefly recall their results. As usual we denote
by $S$ the closure of the operator:
$$S(x\Omega)=x^{*}\Omega\;\;\;{\rm\; for\; all\;} x\in \Gamma_{q}(H_{\R},(U_{t})_{t \in \R})$$
Let $S=J\Delta^{\frac{1}{2}}$ be its polar decomposition. $J$ and
$\Delta$ are respectively the modular conjugation and the modular
operator relative to $\varphi$. The following explicit formulas hold
:
$$S(h_{1}\otimes \dots \otimes h_{n})=h_{n}\otimes \dots \otimes
h_{1}\;\;\;{\rm\; for\; all\;} h_{1},\dots,h_{n} \in H_{\R}$$
$\Delta$ is the closure of the operator
$\mathop{\oplus}\limits_{n=0}^{\infty}(A^{-1})^{\otimes n}$ and
$$J(h_{1}\otimes \dots \otimes h_{n})=A^{-\f{1}{2}}h_{n}\otimes
\dots \otimes A^{-\f{1}{2}}h_{1}\;\;\;{\rm\; for\; all\;}
h_{1},\dots,h_{n} \in H_{\R}\cap {\rm dom}A^{-\f{1}{2}}$$ The
modular group of automorphisms $(\sigma_{t})_{t \in \R}$ on
$\Gamma_{q}(H_{\R},(U_{t})_{t \in \R})$ relative to $\varphi$ is
given by:
$$\sigma_{t}(G(f))=\Delta^{it}G(f)\Delta^{-it}=G(U_{-t}f)\qquad\text{for all}\quad t\in
\R\quad\text{and all}\quad f\in H_{\R}$$ In the following Lemma we
state a well known formula giving, in particular, all moments of the
$q$-Gaussians.
\begin{lem}\label{mom} Let $r \in \N_{*}$ and $(h_{l})_{\substack{-r\le l \le r\\
k\neq 0}}$ be a family of vectors in $H_{\R}$. For all $l\in
\{1,\dots, r\}$ consider the operator
$d_{l}=a^{*}(h_{l})+a(h_{-l})$. For all $(k(1),\dots, k(r))\in
\{1,*\}^{r}$ we have:

\begin{equation*}
\varphi(d_{1}^{k(1)}\dots d_{r}^{k(r)}) = \left\{
\begin{array}{cl} 0&
\text{if\;\;} r \text{\;is odd}\vspace{0.5
cm}\\
\sum\limits_{\substack{
\mathcal{V}\text{\;2-partition}\\
\mathcal{V}=\{(s_{l},t_{l})_{l=1}^{l=p}\} \text{with\;}s_{l}<t_{l}}}
q^{i(\mathcal{V})}\Prod_{l=1}^{p}\varphi(d_{s_{l}}^{k(s_{l})}d_{t_{l}}^{k(t_{l})})
& \text{if\;\;} r=2p
\end{array}
\right.\\
\end{equation*}
 where $i(\mathcal{V})=\#\{(k,l),\;s_{k}<s_{l}<t_{k}<t_{l}\}$ is
the number of crossings of the $2$-partition $\mathcal{V}$.
\end{lem}

\medskip

\noindent{\it Remarks.}

\smallskip

{\noindent $\bullet$ When} $(U_{t})_{t \in \R}$ is trivial,
$\Gamma_{q}(H_{\R},(U_{t})_{t \in \R})$ reduces to Bo${\rm
\dot{z}}$ejko and Speicher's $q$-Gaussian algebra
$\Gamma_{q}(H_{\R})$. This is the only case where $\varphi$ is a
trace on $\Gamma_{q}(H_{\R},(U_{t})_{t \in \R})$. Actually,
$\Gamma_{q}(H_{\R})$ is known to be a non-hyperfinite $I\!I_{1}$
factor (see \cite{BKS}, \cite{BS1}, \cite{Nou} and \cite{Ri}). In
all other cases $\Gamma_{q}(H_{\R},(U_{t})_{t \in \R})$ turns out to
be a type $I\!I\!I$ von Neumann algebra (see \cite{Sh} and
\cite{Hi}).

\smallskip

{\noindent $\bullet$} Lemma \ref{mom} implies that for all $n\in \N$
and all $f\in H_{\R}$:
$$\varphi(G(f)^{2n})=\sum\limits_{
\mathcal{V}\text{\;2-partition}}
q^{i(\mathcal{V})}\|f\|_{H_{\R}}^{2n}$$ Therefore, we see that the
distribution of a single gaussian does not depend on the group
$(U_{t})_{t \in \R}$. In the tracial case (thus in all cases), and
when $\|f\|=1$, this distribution is the absolutely continuous
probability measure $\nu_{q}$ supported on the interval
$[-2/\sqrt{1-q},2/\sqrt{1-q}]$ whose orthogonal polynomials are the
$q$-Hermite polynomials (see \cite{BKS}). In particular, we have :
\begin{equation}\label{moman}\text{For
all}\quad f\in H_{\R},\quad\|G(f)\|=\f{2}{\sqrt{1-q}}\|f\|_{H_{\R}}
\end{equation}

\subsection{The finite dimensional case}\label{subsec2}

We now briefly recall a description of the von Neumann algebra
$\Gamma_{q}(H_{\R},U_{t})$ where $H_{\R}$ is an Euclidian space of
dimension $2k$ ($k \in \N_{*}$). There exists $(H_{j})_{1\le j \le
k}$ a family of two dimensional spaces, invariant under $(U_{t})_{t
\in \R}$, and $(\lambda_{j})_{1\le j \le k}$ some real numbers
greater or equal to $1$ such that for all $ j \in \{1,\dots, k\}$,
$$H_{\R}=\mathop{\oplus}_{1\le j \le k}H_{j}\quad\text{and}
\quad U_{t|H_{j}}= \left(
\begin{array}{cc}
\cos(t\ln(\lambda_{j})) & -\sin(t\ln(\lambda_{j}))\\
\sin(t\ln(\lambda_{j})) & \cos(t\ln(\lambda_{j}))
\end{array}
\right)$$ We put $I=\{-k,\dots,-1\}\cup\{1,\dots,k\}$. It is then
easily checked that the deformed scalar product
$\scal{\,.\,}{\,.\,}_{U}$ on the complexification $H_{\C}$ of
$H_{\R}$ is characterized by the condition that there exists a basis
$(f_{j})_{j \in I}$ in $H_{\R}$ such that for all $(j,l) \in
\{1,\dots,k\}^{2}$

\begin{equation}\label{cond}
\scal{f_{j}}{f_{-l}}_{U}=\delta_{j,l}.i\f{\lambda_{j}-1}{\lambda_{j}+1}
\quad\text{and}\quad\scal{f_{\pm j}}{f_{\pm l}}_{U}=\delta_{j,l}
\end{equation}
 For all $j\in \{1,\dots,k\}$ we put
$\mu_{j}=\lambda_{j}^{\f{1}{4}}$. Let $(e_{j})_{j \in I}$ be a real
orthonormal basis of $\C^{2k}$ equipped with its canonical scalar
product. For all $j \in \{1,\dots,k\}$ we put
$$\hat{f}_{j}=\f{1}{\sqrt{\mu_{j}^{2}+\mu_{j}^{-2}}}\left(\mu_{j}e_{-j}+\mu_{j}^{-1}e_{j}\right)\quad\text{and}\quad
\hat{f}_{-j}=\f{i}{\sqrt{\mu_{j}^{2}+\mu_{j}^{-2}}}\left(\mu_{j}e_{-j}-\mu_{j}^{-1}e_{j}\right)$$
It is easy to see that the conditions (\ref{cond}) are fulfilled for
the family $(\hat{f}_{j})_{j\in I}$. We will denote by $H_{\R}$ the
Euclidian space generated by the family $(\hat{f}_{j})_{j\in I}$ in
$\C^{2k}$. This provides us with a realization of
$\Gamma_{q}(H_{\R},U_{t})$ as a subalgebra of
$B\left(\mathcal{F}_{q}(\C^{2k})\right)$. Indeed,
$\Gamma_{q}(H_{\R},U_{t})=\{G(\hat{f}_{j}),\,j\in I\}^{\prime
\prime}\subset B\left(\mathcal{F}_{q}(\C^{2k})\right)$. For all $j
\in \{1,\dots,k\}$ put
$$f_{j}=\f{\sqrt{\mu_{j}^{2}+\mu_{j}^{-2}}}{2}\hat{f}_{j}\quad\text{and}
\quad f_{-j}=\f{\sqrt{\mu_{j}^{2}+\mu_{j}^{-2}}}{2}\hat{f}_{-j}$$ We
define the following generalized semi-circular variable by:
$$c_{j}=G(f_{j})+iG(f_{-j})=W(f_{j}+if_{-j})$$
It is clear that $\Gamma_{q}(H_{\R},U_{t})=\{c_{j},\;j\in
\{1,\dots,k\}\}''\subset B\left(\mathcal{F}_{q}(\C^{2k})\right)$ and
we can check that
\begin{equation}\label{semicirc}
c_{j}=\mu_{j}a(e_{-j})+\mu_{j}^{-1}a^{*}(e_{j})
\end{equation}
Moreover, for all $j \in \{1,\dots,k\}$, $c_{j}$ is an entire vector
for $(\sigma_{t})_{t \in \R}$ and we have, for all $z \in \C$:
$$\sigma_{z}(c_{j})=\lambda_{j}^{iz}c_{j}.$$
\medskip
Recall that all odd $*$-moments of the family $(c_{j})_{1\le j \le
k}$ are zero. Applying Lemma \ref{mom} to the operators $c_{j}$ we
state, for further references, an explicit formula for the
$*$-moments of $(c_{j})_{1\le j \le k}$. In the following we use the
convention $c^{-1}=c^{*}$ when there is no possible confusion.

\begin{lem}\label{momdef} Let $r \in \N_{*}$, $(j(1),\dots,j(2r))\in
\{1,\dots,k\}^{2r}$ and $(k(1),\dots,k(2r))\in \{\pm 1\}^{2r}$
\begin{align*}
\varphi(c_{j(1)}^{k(1)}\dots c_{j(2r)}^{k(2r)})& =\sum_{\substack{
\mathcal{V}\text{\; 2-partition}\\
\mathcal{V}=\{(s_{l},t_{l})_{l=1}^{l=r}\}
\text{with\;}s_{l}<t_{l}}}q^{i(\mathcal{V})}
\Prod_{l=1}^{r}\varphi(c^{k(s_{l})}_{j(s_{l})}c_{j(t_{l})}^{k(t_{l})})\\
& =\sum_{\substack{ \mathcal{V}\text{\;2-partition}\\
\mathcal{V}=\{(s_{l},t_{l})_{l=1}^{l=r}\}
\text{with\;}s_{l}<t_{l}}}q^{i(\mathcal{V})}
\Prod_{l=1}^{r}\mu_{j(s_{l})}^{2k(s_{l})}\,\delta_{k(s_{l}),\;-k(t_{l})}\,\delta_{j(s_{l}),j(t_{l})}
\end{align*}
\end{lem}

\proof As said above this is a consequence of Lemma \ref{mom} and
the explicit computation of covariances. Using (\ref{semicirc}) we
have:
\begin{align*}
\varphi(c^{k(1)}_{j(1)}c_{j(2)}^{k(2)})&=\scal{c^{-k(1)}_{j(1)}\Omega}
{c_{j(2)}^{k(2)}\Omega}\\
&=\scal{\mu_{j(1)}^{k(1)}e_{-k(1)j(1)}}{\mu_{j(2)}^{-k(2)}e_{k(2)j(2)}}\\
&=\mu_{j(1)}^{2k(1)}\,\delta_{k(1),\;-k(2)}\,\delta_{j(1),j(2)}
\end{align*}
\Eproof

\subsection{Baby Fock}\label{subsec3}

The symmetric Baby Fock (also known as symmetric toy Fock space) is
at some point a discrete approximation of the bosonic Fock space
(see \cite{PAM}). In \cite{Bi}, Biane considered spin systems with
mixed commutation and anti-commutation relations (which is a
generalization of the symmetric toy Fock), and used it to
approximate $q$-Fock space (via Speicher central limit Theorem). In
this section we recall the formal construction of \cite{Bi}. Let $I$
be a finite subset of $\Z$ and $\epsilon$ a function from $I\times
I$ to $\{-1,1\}$ satisfying for all $(i,j)\in I^{2}$,
$\epsilon(i,j)=\epsilon(j,i)$ and $\epsilon(i,i)=-1$. Let
$\mathcal{A}(I,\epsilon)$ be the free complex unital algebra with
generators $(x_{i})_{i \in I}$
 quotiented by the relations
\begin{equation}\label{mix}
x_{i}x_{j}-\epsilon(i,j) x_{j}x_{i}=2\delta_{i,j}\quad{\text
for}\quad (i,j)\in I^{2}
\end{equation}
We define an involution on $\mathcal{A}(I,\epsilon)$ by
$x_{i}^{*}=x_{i}$. For a subset $A=\{i_{1},\dots ,i_{k}\}$ of $I$
with $i_{1}<\dots <i_{k}$ we put $x_{A}=x_{i_{1}} \dots x_{i_{k}}$,
where, by convention, $x_{\emptyset}=1$. Then $(x_{A})_{A\subset I}$
is a basis of the vector space $\mathcal{A}(I,\epsilon)$. Let
$\varphi^{\epsilon}$ be the tracial functional defined by
$\varphi^{\epsilon}(x_{A})=\delta_{A, \emptyset}$ for all $A \subset
I$. $\scal{x}{y}=\varphi^{\epsilon}(x^{*}y)$ defines a positive
definite hermitian form on $\mathcal{A}(I,\epsilon)$. We will denote
by $L^{2}(\mathcal{A}(I,\epsilon),\varphi^{\epsilon})$ the Hilbert
space $\mathcal{A}(I,\epsilon)$ equipped with $\scal{\,.\,}{\,.\,}$.
$(x_{A})_{A\subset I}$ is an orthonormal basis of
$L^{2}(\mathcal{A}(I,\epsilon),\varphi^{\epsilon})$. For each $i \in
I$, define the following partial isometries $\beta_{i}^{*}$ and
$\alpha_{i}^{*}$ of
$L^{2}(\mathcal{A}(I,\epsilon),\varphi^{\epsilon})$ by:
$$\beta_{i}^{*}(x_{A})=\left\{
\begin{array}{cl}
x_{i}x_{A} & {\rm if}\quad i \not\in A\\
0 & {\rm if}\quad i\in A
\end{array}
\right. \quad\text{and}\quad\alpha_{i}^{*}(x_{A})=\left\{
\begin{array}{cl}
x_{A}x_{i} & {\rm if}\quad i \not\in A\\
0 & {\rm if}\quad i\in A
\end{array}
\right.
$$
Note that their adjoints are given by:
$$\beta_{i}(x_{A})=\left\{
\begin{array}{cl}
x_{i}x_{A} & {\rm if}\quad i\in A\\
0 & {\rm if}\quad i\not\in A
\end{array}
\right. \quad\text{and}\quad\alpha_{i}(x_{A})=\left\{
\begin{array}{cl}
x_{A}x_{i} &{\rm if}\quad i\in A\\
0 & {\rm if}\quad i\not\in A
\end{array}
\right.
$$
$\beta_{i}^{*}$ and $\beta_{i}$ (respectively $\alpha_{i}^{*}$ and
$\alpha_{i}$) are called the left (respectively right) creation and
annihilation operators at the Baby Fock level. In the next Lemma we
recall from \cite{Bi} the fundamental relations 1. and 2., and we
leave the proof of 3., 4. and 5. to the reader.

\begin{lem}\label{relcrea}The following relations hold:
\begin{enumerate}
\item For all $i\in I$ $(\beta_{i}^{*})^{2}=\beta_{i}^{2}=0$  and
$\beta_{i}\beta_{i}^{*}+\beta_{i}^{*}\beta_{i}=Id$.
\item For all $(i,j)\in I^{2}$ with $i\neq j$
$\beta_{i}\beta_{j}-\epsilon(i,j)\beta_{j}\beta_{i}=0$ and
$\beta_{i}\beta_{j}^{*}-\epsilon(i,j)\beta_{j}^{*}\beta_{i}=0$.
\item Same relations as in 1. and 2. with $\alpha$ in place of $\beta$.
\item For all $i\in I$
$\beta_{i}^{*}\alpha_{i}^{*}=\alpha_{i}^{*}\beta_{i}^{*}=0$ and for
all $(i,j)\in I^{2}$ with $i\neq j$
$\beta_{i}^{*}\alpha_{j}^{*}=\alpha_{j}^{*}\beta_{i}^{*}$.
\item For all $(i,j)\in I^{2}$
$\beta_{i}^{*}\alpha_{j}=\alpha_{j}\beta_{i}^{*}$.
\end{enumerate}
\end{lem}

It is easily seen, by 1. and 2. of Lemma \ref{relcrea}, that the
self adjoint operators defined by:
$\gamma_{i}=\beta_{i}^{*}+\beta_{i}$ satisfy the following relation
:
\begin{equation}\label{eq13}\text{for all}\quad (i,j)\in
I^{2},\quad\gamma_{i}\gamma_{j}-\epsilon(i,j)\gamma_{j}\gamma_{i}=2\delta_{i,j}\text{Id}
\end{equation}

Let $\Gamma_{I}\subset
B(L^{2}(\mathcal{A}(I,\epsilon),\varphi^{\epsilon}))$ be the
$*$-algebra generated by all $\gamma_{i}$, $i\in I$. Still denoting
by $\varphi^{\epsilon}$ the vector state associated to the vector
$1$, it is known that $\varphi^{\epsilon}$ is a faithful normalized
trace on the finite dimensional $C^{*}$-algebra $\Gamma_{I}$ (see
the remarks below). Moreover, $\Gamma_{I}\subset
B(L^{2}(\mathcal{A}(I,\epsilon),\varphi^{\epsilon}))$ is the
faithful GNS representation of $(\Gamma_{I},\varphi^{\epsilon})$
with cyclic and separating  vector $1$.

\bigskip

\noindent{{\it Remarks.}}

\smallskip

{\noindent $\bullet$ It is clear} that we can do the previous
construction for some finite sets $I$ that are not given explicitly
as subsets of $\Z$. Then, to each total order on $I$ we can
associate a basis $(x_{A})_{A\subset I}$ of
$\mathcal{A}(I,\epsilon)$. But, because of the commutation relations
(\ref{mix}), the state $\varphi^{\epsilon}$, the scalar product on
$\mathcal{A}(I,\epsilon)$ and the creation operators do not depend
on the chosen total order.

\smallskip

{\noindent $\bullet$ We can} also extend the previous construction
to not necessarily finite sets $I$. Only the faithfulness of
$\varphi^{\epsilon}$ on $\Gamma_{I}$ requires some comments. It
suffices to see that the vector $1$ is separating for $\Gamma_{I}$.
Indeed, set $\delta_{i}=\alpha_{i}^{*}+\alpha_{i}$ for $i \in I$,
and
$$\Gamma_{r,I}=\{\delta_{i}, i\in I\}^{\prime \prime}\subset B(L^{2}(\mathcal{A}(I,\epsilon),\varphi^{\epsilon})).$$
Then, it is clear from 4. and 5. of Lemma \ref{relcrea}, that
$\Gamma_{r,I}\subset \Gamma_{I}^{\prime}$ (there is actually
equality). Since $1$ is clearly cyclic for $\Gamma_{r,I}$, then it
is also cyclic for $\Gamma_{I}^{\prime}$, thus $1$ is separating for
$\Gamma_{I}$.

\smallskip

{\noindent $\bullet$ Let $I$} and $J$, $I\subset J$, be some sets
together with signs $\epsilon$ and $\epsilon^{\prime}$ such that
$\epsilon^{\prime}_{|I\times I}=\epsilon$. It is clear that
$L^{2}(\mathcal{A}(I,\epsilon),\varphi^{\epsilon})$ embeds
isometrically in
$L^{2}(\mathcal{A}(J,\epsilon^{\prime}),\varphi^{\epsilon^{\prime}})$.
Set $K=J\smallsetminus I$. Fix some total orders on $I$ and $K$ and
consider the total order on $J$ which coincides with the orders of
$I$ and $K$ and such that any element of $I$ is smaller than any
element of $K$. The associated orthonormal basis of
$L^{2}(\mathcal{A}(J,\epsilon^{\prime}),\varphi^{\epsilon^{\prime}})$
is given by the family $(x_{A}x_{B})_{A\in \mathcal{F}(I), B\in
\mathcal{F}(K)}$ (where $\mathcal{F}(I)$, respectively
$\mathcal{F}(K)$, denotes the set of finite subsets of $I$,
respectively $K$). In particular with have the following Hilbertian
decomposition:
\begin{equation}\label{deco}
L^{2}(\mathcal{A}(J,\epsilon^{\prime}),\varphi^{\epsilon^{\prime}})
=\bigoplus_{B\in
\mathcal{F}(K)}L^{2}(\mathcal{A}(I,\epsilon),\varphi^{\epsilon})x_{B}
\end{equation} For $j\in I$ we (temporarily) denote by
$\widetilde{\beta}_{j}$ the annihilation operator in
$B(L^{2}(\mathcal{A}(I,\epsilon),\varphi^{\epsilon}))$ and simply by
$\beta_{j}$ its analogue in
$B(L^{2}(\mathcal{A}(J,\epsilon^{\prime}),\varphi^{\epsilon^{\prime}}))$.
Let $\widetilde{C}_{I}$ (respectively $C_{J}$) be the
$C^{*}-$algebra generated by $\{\widetilde{\beta}_{j},\;j\in I\}$
(respectively $\{\beta_{j},\;j\in J\}$) in
$B(L^{2}(\mathcal{A}(I,\epsilon),\varphi^{\epsilon}))$ (respectively
$B(L^{2}(\mathcal{A}(J,\epsilon^{\prime}),\varphi^{\epsilon^{\prime}}))$).
Consider also $C_{I}$ the $C^{*}-$algebra generated by
$\{\beta_{j},\;j\in I\}$ in \\
$B(L^{2}(\mathcal{A}(I,\epsilon),\varphi^{\epsilon}))$. For
$B=\{j_{1},\dots,j_{k}\}\subset K$, with $j_{1}<\dots<j_{k}$, let us
denote by $\alpha_{B}$ the operator
$\alpha_{j_{1}}\dots\alpha_{j_{k}}$. If $\widetilde{T}\in
\widetilde{C}_{I}$ and if $T$ denotes its counterpart in $C_{I}$,
then it is easily seen that, with respect to the Hilbertian
decomposition (\ref{deco}), we have
\begin{equation}\label{comp}
T=\bigoplus_{B\in
\mathcal{F}(K)}\alpha_{B}^{*}\widetilde{T}\alpha_{B}.
\end{equation}
It follows that $\widetilde{C}_{I}$ is $*$-isomorphic to
$C_{I}\subset C_{J}$.

\smallskip

{\noindent $\bullet$ It is possible} to find explicitly selfadjoint
matrices satisfying the mixed commutation and anti-commutation
relations (\ref{eq13}) (see \cite{Sp} and \cite{Bi}). We choose to
present this approach because it will be easier to handle the
objects of modular theory in this abstract situation when we will
deal with non-tracial von Neumann algebras (see section
\ref{section4}).

\subsection{Speicher's central limit Theorem}\label{subsec4}

We recall Speicher's central limit theorem which is specially
designed to handle either commuting or anti-commuting (depending on
a function $\epsilon$) independent variables. Roughly speaking,
Speicher's central limit theorem asserts that such a family of
centered noncommutative variables which have a fixed covariance, and
uniformly bounded $*$-moments, is convergent in $*$-moments, as soon
as a combinatorial quantity associated with $\epsilon$ is
converging. Moreover the limit $*$-distribution is only determined
by the common covariance and the limit of the combinatorial
quantity.\\
We start by recalling some basic notions on independence and set
partitions.
\begin{defi} Let $(\mathcal{A},\varphi)$ be a $*$-algebra equipped with a state $\varphi$
and $(\mathcal{A}_{i})_{i\in I}$ a family of $C^{*}$-subalgebras of
$\mathcal{A}$. The family $(\mathcal{A}_{i})_{i\in I}$ is said to be
independent if for all $r \in \N_{*}$, $(i_{1},\dots,i_{r})\in
I^{r}$ with $i_{s}\neq i_{t}$ for $s\neq t$, and all $a_{i_{s}} \in
\mathcal{A}_{i_{s}}$ for $s \in \{1,\dots,r\}$ we have:
$$\varphi(a_{i_{1}}\dots a_{i_{r}})=\varphi(a_{i_{1}})\dots
\varphi(a_{i_{r}})$$ As usual, a family $(a_{i})_{i \in I}$ of
non-commutative random variables of $\mathcal{A}$ will be called
independent if the family of $C^{*}$-subalgebras of $\mathcal{A}$
that they generate is independent.
\end{defi}
On the set of $p$-uples of integers belonging to $\{1,\dots,N\}$
define the equivalence relation $\sim$ by:
$$(i(1),\dots, i(p))\sim (j(1),\dots, j(p))\;\;\text{if}\;\;
\big(i(l)=i(m) \Longleftrightarrow j(l)=j(m)\big)\;\forall\;(l,m)\in
\{1,\dots,p\}^{2}$$ Then the equivalence classes for the relation
$\sim$ are given by the partitions of the set $\{1,\dots,p\}$. We
denote by $V_{1},\dots,V_{r}$ the blocks of the partition
$\mathcal{V}$ and we call $\mathcal{V}$ a $2$-partition if each of
these blocks is of cardinal $2$.  The set of all $2$-partitions of
the set $\{1,\dots,p\}$ ($p$ even) will be denoted by
$\mathcal{P}_{2}(1,\dots,p)$ . For $\mathcal{V}\in
\mathcal{P}_{2}(1,\dots,2r)$  let us denote by
$V_{l}=(s_{l},t_{l})$, $s_{l}< t_{l}$, for $l\in \{1,\dots,r\}$ the
blocks of the partition $\mathcal{V}$. The set of crossings of
$\mathcal{V}$ is defined by
$$I(\mathcal{V})=\{(l,m)\in \{1,\dots,r\}^{2},\;
s_{l}<s_{m}<t_{l}<t_{m}\}$$ The $2$-partition $\mathcal{V}$ is said
to be crossing if $I(\mathcal{V})\neq \emptyset$ and non-crossing if
$I(\mathcal{V})= \emptyset$.
\begin{thm}[Speicher]\label{TCL}
Consider $k$ sequences $(b_{i,j})_{(i,j)\in
\N_{*}\times\{1,\dots,k\}}$ in a non-\\commutative probability space
$(B,\varphi)$ satisfying the following conditions:
\begin{enumerate}
\item The family $(b_{i,j})_{(i,j)\in
\N_{*}\times\{1,\dots,k\}}$ is independent.
\item For all $(i,j)\in
\N_{*}\times\{1,\dots,k\}$, $\varphi(b_{i,j})=0$
\item For all $(k(1),k(2)) \in \{-1,1\}^{2}$ and $(j(1),j(2))\in \{1,\dots,k\}^{2}$, the covariance\\
$\varphi(b_{i,j(1)}^{k(1)}b_{i,j(2)}^{k(2)})$ is independent of $i$
and will be denoted by $\varphi(b_{j(1)}^{k(1)}b_{j(2)}^{k(2)})$.
\item For all $w \in \N_{*}$, $(k(1),\dots,k(w))\in \{-1,1\}^{w}$
and all $j \in \{1,\dots,k\}$ there exists a constant $C$ such that
for all $i\in \N_{*}$, $|\varphi(b_{i,j}^{k(1)}\dots
b_{i,j}^{k(w)})|\le C$.
\item For all $(i(1),i(2)) \in \N_{*}^{2}$ there exists a sign
$\epsilon(i(1),i(2)) \in \{-1,1\}$ such that for all $(j(1),j(2))
\in \{1,\dots,k\}^{2}$ with $(i(1),j(1))\neq (i(2),j(2))$ and all
$(k(1),k(2))\in \{-1,1\}^{2}$ we have
$$b_{i(1),j(1)}^{k(1)}b_{i(2),j(2)}^{k(2)}-\epsilon(i(1),i(2))b_{i(2),j(2)}^{k(2)}b_{i(1),j(1)}^{k(1)}=0.$$
(notice that the function $\epsilon$ is necessarily symmetric in its
two arguments).
\item For all $r \in \N_{*}$ and all $\mathcal{V}=\{(s_{l},t_{l})_{l=1}^{l=r}\} \in
\mathcal{P}_{2}(1,\dots,2r)$  the following limit exists
$$t(\mathcal{V})=\lim_{N\rightarrow
+\infty}\f{1}{N^{r}}\sum_{\substack{i(s_{1}),\dots, i(s_{r})=1\\
i(s_{l})\neq i(s_{m})\text{\;for\;} l\neq m}}^{N}\Prod_{(l,m)\in
I(\mathcal{V})}\epsilon\left(i(s_{l}),i(s_{m}))\right)$$

Let $S_{N,j}=\f{1}{\sqrt{N}}\sum_{i=1}^{N}b_{i,j}$. Then we have for
all $p \in \N_{*}$, $(k(1),\dots,k(p))\in \{-1,1\}^{p}$ and all
$(j(1),\dots,j(p)) \in \{1,\dots,k\}^{p}$:
$$
\lim_{N\rightarrow +\infty}\varphi(S_{N, j(1)}^{k(1)}\dots S_{N,
j(p)}^{k(p)})=\left\{
\begin{array}{cl}
0 &\text{if $p$ is odd}\vspace{0.5 cm}\\

\displaystyle{\sum_{\substack{\mathcal{V} \in \mathcal{P}_{2}(1,\dots,2r)\\
\mathcal{V}=\{(s_{l},t_{l})_{l=1}^{l=r}\}
}}}t(\mathcal{V})\Prod_{l=1}^{r}
\varphi(b_{j(s_{l})}^{k(s_{l})}b_{j(t_{l})}^{k(t_{l})}) &
\text{if\;} p=2r
\end{array}
\right.$$
\end{enumerate}
\end{thm}

\noindent{\it Remark.} Speicher's Theorem is proved in \cite{Sp} for
a single limit variable. One could either convince oneself that the
proof of Theorem \ref{TCL} goes along the same lines, or deduce it
from Speicher's usual theorem. Indeed, it suffices to apply
Speicher's theorem to the family
$\left(\sum\limits_{j=1}^{k}z_{j}b_{i,j}\right)_{i \in \N}$, for all
$(z_{1},\dots,z_{k})\in \T^{k}$ and to identify the Fourier
coefficients of the limit $*$-moments.

The following Lemma, proved in \cite{Sp}, guarantees the almost sure
convergence of the quantity $t(\mathcal{V})$ provided that the
function $\epsilon$ has independent entries following the same
$2$-points Dirac distribution:

\begin{lem}\label{rand} Let $q \in (-1,1)$ and consider a family of random variables
$\epsilon(i,j)$ for $(i,j)\in \N_{*}$ with $i\neq j$, such that
\begin{enumerate}
\item For all $(i,j)\in \N_{*}$ with $i\neq
j$, $\epsilon(i,j)=\epsilon(j,i)$
\item The family $(\epsilon(i,j))_{i>j}$ is independent
\item For all $(i,j)\in
\N_{*}$ with $i\neq j$ the probability distribution of
$\epsilon(i,j)$ is
$$\f{1+q}{2}\delta_{1}+\f{1-q}{2}\delta_{-1}$$
\end{enumerate}
Then, almost surely, we have for all $r\in \N_{*}$ and for all
$\mathcal{V} \in \mathcal{P}_{2}(1,\dots,2r)$
$$\lim_{N\rightarrow
+\infty}\f{1}{N^{r}}\sum_{\substack{i(s_{1}),\dots, i(s_{r})=1\\
i(s_{l})\neq i(s_{m})\text{\;for\;} l\neq m}}^{N}\Prod_{(l,m)\in
I(\mathcal{V})}\epsilon(i(s_{l}),i(s_{m}))=q^{i(\mathcal{V})}$$
\end{lem}

\noindent{\it Remark.} It is now a straightforward verification to
see that Theorem \ref{TCL} combined with Lemma \ref{rand} can be
applied to families of mixed commuting /anti-commuting Gaussian
operators (see Lemma \ref{ind} for the independence condition). The
limit moments are those given by the  classical $q$-Gaussian
operators (by classical we mean that $(U_{t})_{t \in \R}$ is
trivial).\\
Alternatively, one can apply directly Speicher's theorem  to
families of mixed commuting /anti-commuting creation operators as it
is done in \cite{Sp} and \cite{Bi}. The limit $*$-moments are in
this case the $*$-moments of classical $q$-creation operators.

\section{The tracial case}\label{section2}

Our goal in this section is to show that $\Gamma_{q}(H_{\R})$ is
QWEP. In fact, by inductive limit, it is sufficient to prove it for
$H_{\R}$ finite dimensional. Let $k \ge 1$. We will consider
$\R^{k}$ as the real Hilbert space of dimension $k$, with the
canonical orthonormal basis $(e_{1}, \dots , e_{k})$, and $\C^{k}$,
its complex counterpart. Let us fix $q \in (-1,1)$ and consider
$\Gamma_{q}(\R^{k})$ the von Neumann algebra generated by the
$q$-Gaussians $G(e_{1}), \dots, G(e_{k})$. We denote by $\tau$ the
expectation on the vacuum vector, which is a trace in this
particular case.

By the ending remark of section \ref{section1}., there are Hermitian
matrices, $g_{n, 1}(\omega), \dots , g_{n, k}(\omega)$, depending on
a random parameter denoted by $\omega$ and lying in a finite
dimensional matrix algebra, such that their joint $*$-distribution
converges almost surely to the joint $*$-distribution of the
$q$-Gaussians in the following sense: for all polynomial $P$ in $k$
noncommuting variables,

$$\lim_{n \to \infty} \tau_{n}(P(g_{n, 1}(\omega), \dots , g_{n, k}(\omega)))=
\tau (P(G(e_{1}), \dots , G(e_{k}))) \hspace{0.5 cm} {\rm almost\;
surely\;in\; \omega}.$$ We will denote by $\mathcal{A}_{n}$ the
finite dimensional $C^*$-algebra generated by $g_{n, 1}(\omega),
\dots ,g_{n, k}(\omega)$. We recall that these algebras are equipped
with the trace $\tau_{n}$ defined by:
$$\tau_{n}(x)=\scal{1}{x.1}$$
Since the set of all monomials in $k$ noncommuting variables is
countable, we have for almost all $\omega$,
\begin{equation}\label{eq0}
\lim_{n \to \infty} \tau_{n}(P(g_{n, 1}(\omega), \dots , g_{n,
k}(\omega)) =\tau (P(G(e_{1}), \dots , G(e_{k}))) \hspace{0.5 cm}
{\rm for\; all\; such\; monomials\; }P
\end{equation}
A fortiori we can find an $\omega_{0}$ such that (\ref{eq0}) holds
for $\omega_{0}$.
 We will fix such an $\omega_{0}$ and simply denote by $g_{n, i}$ the matrix
 $g_{n, i}(\omega_{0})$ for all $i \in \{1,\dots , k\}$. With these notations,
 it is clear that, by linearity, we have for all polynomials $P$ in $k$ noncommuting variables,

\begin{equation}\label{eq}
\lim_{n \to \infty} \tau_{n}(P(g_{n, 1}, \dots , g_{n, k}))=\tau
(P(G(e_{1}), \dots , G(e_{k}))).
\end{equation}
We need to have a uniform control on the norms of the matrices
$g_{n, i}$. Let $C$ be such that $\|G(e_{1})\| < C$, we will replace
the $g_{n, i}$'s by their truncations $\chi_{]-C,C[}(g_{n, i})g_{n,
i}$ (where $\chi_{]-C,C[}$ denotes the characteristic function of
the interval  $]-C,C[$). For simplicity $\chi_{]-C,C[}(g_{n,
i})g_{n, i}$ will be denoted by $\tilde{g}_{n, i}$. We now check
that (\ref{eq}) is still valid for the $\tilde{g}_{n, i}$'s.

\begin{lem}\label{pres} With the notations above, for all polynomials $P$ in $k$ noncommuting variables we have

\begin{equation}\label{eq2}
\lim_{n \to \infty} \tau_{n}(P(\tilde{g}_{n, 1}, \dots ,
\tilde{g}_{n, k}))=\tau (P(G(e_{1}), \dots , G(e_{k}))).
\end{equation}
\end{lem}

\proof We just have to prove that for all monomials $P$ in $k$
noncommuting variables we have
$$\lim_{n \to \infty} \tau_{n}\left[P(\tilde{g}_{n, 1}, \dots , \tilde{g}_{n, k})
-P(g_{n, 1}, \dots , g_{n, k})\right]=0.$$ Writing $g_{n,
i}=\tilde{g}_{n, i}+(g_{n, i}-\tilde{g}_{n, i})$ and developing
using multilinearity, we are reduced to showing that the
$L^{1}$-norms of any monomial in $\tilde{g}_{n, i}$ and $(g_{n,
i}-\tilde{g}_{n, i})$ (with at least one factor $(g_{n,
i}-\tilde{g}_{n, i})$) tend to $0$. By the H\"older inequality and
the uniform boundedness of the $\|\tilde{g}_{n, i}\|$'s, it suffices
to show that for all $i \in \{1 \dots k\}$,
\begin{equation}\label{eq3}
\lim_{n \to \infty} \tau_{n}(|\tilde{g}_{n, i}-g_{n,
i}|^{p})=0\hspace{0.5 cm} {\rm for\;all\;}p\ge 1.
\end{equation}
Let us prove (\ref{eq3}) for $i=1$. We are now in a commutative
setting. Indeed, let us introduce the spectral resolutions of
identity, $E^{n}_{t}$ (respectively $E_{t}$), of $g_{n, 1}$
(respectively $G(e_{1})$). By (\ref{eq}) we have for all polynomials
$P$
$$\lim_{n \to \infty} \tau_{n}(P(g_{n, 1}))=\tau(P(G(e_{1}))).$$
We can rewrite this as follows: for all polynomials $P$
$$\lim_{n \to \infty}\int_{\sigma(g_{n, 1})}P(t)\d\scal{E^{n}_{t}
.1}{1}=\int_{\sigma(G(e_{1}))}P(t)\d\scal{E_{t}.\Omega}{\Omega}.$$
Let $\mu_{n}$ (respectively $\mu$) denote the compactly supported
probability measure $\scal{E^{n}_{t}.1}{1}$ (respectively
$\scal{E_{t}.\Omega}{\Omega}$) on $\R$. With these notations our
assumption becomes: for all polynomials $P$
\begin{equation}\label{eq4}
\lim_{n \to \infty}\int P\d \mu_{n}=\int P\d \mu.
\end{equation}
and (\ref{eq3}) is equivalent to:
\begin{equation}\label{eq5}
\lim_{n \to \infty}\int_{|t|\ge C}|t|^{p}\d \mu_{n}=0\hspace{0.5 cm}
{\rm for\;all\;}p\ge 1.
\end{equation}
Then the result follows from the following elementary Lemma. We give
a proof for sake of completeness.

\begin{lem}\label{conc}Let $(\mu_n)_{n \ge 1} $ be a sequence of compactly supported probability measures
on $\R$ converging in moments to a compactly supported probability
measure $\mu$ on $\R$. Assume that the support of $\mu$ is included
in the open interval $]-C,C[$. Then,
$$\lim_{n \to \infty}\int_{|t|\ge C}\d \mu_{n}=0.$$
Moreover, let $f$ be a borelian function on $\R$ such that there
exist $M>0$ and $r \in \N$ satisfying $|f(t)|\le M(t^{2r}+1)$ for
all $t \ge C$. Then,
$$\lim_{n \to \infty}\int_{|t|\ge C} f \d \mu_{n}=0.$$
\end{lem}

\proof For the first assertion, let $C^{\prime} <C$ such that the
support of $\mu$ is included in $]-C^{\prime},C^{\prime}[$. Let
$\epsilon >0$ and an integer $k$ such that
$\left(\frac{C^{\prime}}{C}\right)^{2k} \le \epsilon$. Let
$P(t)=\left(\frac{t}{C}\right)^{2k}$. It is clear that
$\chi_{\{|t|\ge C\}}(t)\le P(t)$ for all $t \in \R$ and that
$\sup_{|t|<C^{\prime}}P(t)\le \epsilon$. Thus,
$$0\le \limsup_{n\to \infty}\int_{|t|\ge C}\d \mu_{n}\le \lim_{n \to \infty}\int P(t)\d \mu_{n}=
\int P(t)\d \mu\le \epsilon .$$
Since $\epsilon$ is arbitrary, we get $\lim\limits_{n\to \infty}\int_{|t|\ge C}\d \mu_{n}=0$.\\
The second assertion is a consequence of the first one. Let $f$ be a
borelian function on $\R$ such that there exist $M>0$ and $r \in \N$
satisfying $|f(t)|\le M(t^{2r}+1)$ for all $t \in \R$. Using the
Cauchy-Schwarz inequality we get:
\begin{eqnarray*}
 0\le \limsup_{n\to \infty}\int_{|t|\ge C}|f| \d \mu_{n} & \le & \limsup_{n\to \infty}\int_{|t|\ge C}
 M(t^{2r}+1) \d \mu_{n}\\
& \le & M\lim_{n\to \infty}\left(\int(t^{2r}+1)^{2} \d
\mu_{n}\right)^{\frac{1}{2}}\lim_{n \to \infty}
\left(\int_{|t|\ge C}\d \mu_{n}\right)^{\frac{1}{2}}\\
& \le & M\left(\int(t^{2r}+1)^{2} \d \mu\right)^{\frac{1}{2}}\lim_{n
\to \infty} \left(\int_{|t|\ge C}\d \mu_{n}\right)^{\frac{1}{2}}=0
\end{eqnarray*} \Eproof

\medskip

\noindent{\it Remark.} Let us define $\mathcal{A}$ as the
$*$-algebra generated by $G(e_{1}), \dots , G(e_{k})$. Observe that
$\mathcal{A}$ is isomorphic to the $*$-algebra of all polynomials in
$k$ noncommuting variables (the free complex $*$-algebra with $k$
generators). Indeed, if $P(G(e_{1}), \dots , G(e_{k}))=0$ for a
polynomial $P$ in $k$ noncommuting variables, then the equation
$P(G(e_{1}), \dots , G(e_{k}))\Omega=0$ implies that all
coefficients of monomials of highest degree are $0$, and thus $P=0$
by induction. More generally, this remains true for the $q$-Araki
Woods algebras $\Gamma_{q}\left(H_{\R}, (U_{t})_{t \in \R}\right)$:
if $(e_{i})_{i \in I}$ is a free family of vectors in $H_{\R}$ then
the $*$-algebra $\mathcal{A}_{I}$ generated by the family
$(G(e_{i}))_{i \in I}$ is isomorphic to the free complex $*$-algebra
with $I$ generators.

\medskip

Let $\mathcal{U}$ be a free ultrafilter on $\N^{*}$ and consider the
ultraproduct von Neumann algebra (see \cite{P} section 9.10) $N$
defined by
$$N=\left.\left(\prod_{n\ge 1}\mathcal{A}_{n}\right)\right/ I_{\mathcal{U}}$$
where $I_{\mathcal{U}}=\{(x_{n})_{n \ge 1}\in \prod\limits_{n\ge
1}\mathcal{A}_{n},\;\lim \limits_{\mathcal{U}}\tau_{n}(x^*_{n}x_n
)=0\}$. The von Neumann algebra $N$ is equipped with the faithful
normal and normalized trace $\tau((x_{n})_{n\ge
1})=\lim\limits_{\mathcal{U}}\tau_{n}(x_{n})$ (which is well
defined).

Using the asymptotic matrix model for the $q$-Gaussians and by the
preceding remark, we can define a $*$-homomorphism $\varphi$ between
the $*$-algebras $\mathcal{A}$ and $N$ in the following way:
$$\varphi(P(G(e_{1}), \dots , G(e_{k})))=(P(\tilde{g}_{n, 1}, \dots , \tilde{g}_{n, k}))_{n\ge 1} $$
for every polynomial $P$ in $k$ noncommuting variables. By Lemma
\ref{pres}, $\varphi$ is trace preserving on $\mathcal{A}$. Since
the $*$-algebra $\mathcal{A}$ is weak-$*$ dense in
$\Gamma_{q}(\R^{k})$, $\varphi$ extends naturally to a trace
preserving homomorphism of von Neumann algebras, that is still
denoted by $\varphi$ (see Lemma \ref{ext} below for a more general
result). It follows that $\Gamma_{q}(\R^{k})$ is isomorphic to a
sub-algebra of $N$ which is the image of a conditional expectation
(this is automatic in the tracial case). Since the
$\mathcal{A}_{n}$'s are finite dimensional, they are injective,
hence their product is injective and a fortiori has the WEP, and
thus $N$ is QWEP. Since $\Gamma_{q}(\R^{k})$ is isomorphic
 to a sub-algebra of $N$ which is the image of a conditional expectation, $\Gamma_{q}(\R^{k})$
is also QWEP (see \cite{Oz}). We have obtained the following:

\begin{thm}\label{qqwep}
Let $H_{\R}$ be a real Hilbert space and $q \in (-1,1)$. The von
Neumann algebra $\Gamma_{q}(H_{\R})$ is QWEP.
\end{thm}

\proof Our previous discussion implies the result for every finite
dimensional $H_{\R}$. The general result is a consequence of the
stability of QWEP by inductive limit (see \cite{Kir} and \cite{Oz}
Proposition 4.1 (iii)). \Eproof

\medskip

Let $C^{*}_{q}(H_{\R})$ be the $C^{*}$-algebra generated by all
$q$-Gaussians:
$$C^{*}_{q}(H_{\R})=C^{*}\big(\{G(f),f\in H_{\R}\}\big)\subset B(\FQ).$$
We now deduce the following strengthening of Theorem \ref{qqwep}.

\begin{cor}\label{*qwep}
Let $H_{\R}$ be a real Hilbert space and $q \in (-1,1)$. The
$C^{*}$-algebra $C_{q}^{*}(H_{\R})$ is QWEP.
\end{cor}

Let $A$, $B$, with $A\subset B$ be $C^{*}$-algebras. Recall (from
\cite{Oz}) that $A$ is said to be weakly cp complemented in $B$, if
there exists a unital completely positive map $\Phi:
B\longrightarrow A^{**}$ such that $\Phi_{|A}=\id_{A}$. Corollary
\ref{*qwep} is then a consequence of the following Lemma.

\begin{lem}\label{wcp}
The $C^{*}$-algebra $C_{q}^{*}(H_{\R})$ is weakly cp complemented in
the von Neumann algebra $\Gamma_{q}(H_{\R})$.
\end{lem}

\proof For any $t\in \R_{+}$ denote by $\Phi_{t}$ the unital
completely positive maps which are the second quantization of
$e^{-t}\id:H_{\R} \longrightarrow H_{\R}$ (see \cite{BKS}):
$$\Phi_{t}=\Gamma_{q}(e^{-t}\id): \Gamma_{q}(H_{\R})\longrightarrow \Gamma_{q}(H_{\R}),\quad \text{for all}\quad t\ge 0.$$
$(\Phi_{t})_{t\in \R_{+}}$ is a semi-group of unital completely
positive maps which is also known as the $q$-Ornstein-Uhlenbeck
semi-group. By the well-known ultracontractivity of the semi-group
$(\Phi_{t})_{t\in \R_{+}}$ (see \cite{B}),
 for all $t\in \R_{+}^{*}$ and all $W(\xi)\in \Gamma_{q}(H_{\R})$, we have
\begin{equation}\label{ultra}
\|\Phi_{t}(W(\xi))\|\le C_{|q|}^{\frac{3}{2}}
\frac{1}{1-e^{-t}}\|\xi\|.
\end{equation}
On the other hand, as a consequence of the Haagerup-Bo\.zejko's
inequality (see \cite{B}), for every $n\in \N$ and for every
$\xi_{n} \in \HC^{\otimes n}$, we have $W(\xi_{n})\in
C_{q}^{*}(H_{\R})$. Fix $t\in \R_{+}^{*}$, $W(\xi) \in
\Gamma_{q}(H_{\R})$, and write $\xi=\sum\limits_{n\in \N}\xi_{n}$
with $\xi_{n} \in \HC^{\otimes n}$ for all $n$. From our last
observation, for all $N\in \N$,
$$T_{N}=\Phi_{t}\big(W\big(\sum_{n=0}^{N}\xi_{n}\big)\big)=\sum_{n=0}^{N}e^{-tn}W(\xi_{n})\in
C_{q}^{*}(H_{\R}).$$ By (\ref{ultra}), $\Phi_{t}(W(\xi))$ is the
norm limit of the sequence $(T_{N})_{N\in \N}$, so
$\Phi_{t}(W(\xi))$ belongs to $C_{q}^{*}(H_{\R})$. It follows that
$\Phi_{t}$ maps $\Gamma_{q}(H_{\R})$ into $C_{q}^{*}(H_{\R})$.
Moreover, it is clear that
\begin{equation}\label{ida}
\lim_{t\to 0} \|\Phi_{t}(W(\xi))-W(\xi)\|=0,\quad\text{for all}\quad
W(\xi)\in C_{q}^{*}(H_{\R}).
\end{equation}
Take $(t_{n})_{n\in \N}$ a sequence of positive real numbers
converging to $0$ and fix $\mathcal{U}$ a free ultrafilter on $\N$.
By $w^{*}$-compactness of the closed balls in
$\big(C_{q}^{*}(H_{\R})\big)^{**}$, we can define the following
mapping $\Phi:\Gamma_{q}(H_{\R})\longrightarrow
\big(C_{q}^{*}(H_{\R})\big)^{**}$ by
$$\Phi(W(\xi))=w^{*}\text{-}\lim_{n,\mathcal{U}}\Phi_{t_{n}}(W(\xi)), \quad\text{for all}\quad W(\xi)\in \Gamma_{q}(H_{\R}).$$
$\Phi$ is a unital completely positive map satisfying
$\Phi_{|C_{q}^{*}(H_{\R})}=\id_{C_{q}^{*}(H_{\R})}$ by (\ref{ida}).
\Eproof

\bigskip

{\noindent \it Proof of Corollary \ref{*qwep}.} This is a
consequence of Theorem \ref{qqwep}, Lemma \ref{wcp} and Proposition
4.1 (ii) in \cite{Oz}.\Eproof

\section{Embedding into an ultraproduct}\label{section3}

The general setting is as follows. We start with a family
$((\mathcal{A}_{n},\varphi_{n}))_{n \in \N}$ of von Neumann algebras
equipped with normal faithful state $\varphi_{n}$. We assume that
$\mathcal{A}_{n}\subset B(H_{n})$, where the inclusion is given by
the G.N.S. representation of $(\mathcal{A}_{n},\varphi_{n})$. Let
$\mathcal{U}$ be a free ultrafilter on $\N$, and let
$$\widetilde{\mathcal{A}}=\Prod_{n \in
\N}\mathcal{A}_{n}/\mathcal{U}$$ be the $C^{*}$-ultraproduct over
$\mathcal{U}$ of the algebras $\mathcal{A}_{n}$. We canonically
identify $\widetilde{\mathcal{A}}\subset B(H)$,
 where $H=\Prod_{n\in \N}H_{n}/\mathcal{U}$ is the ultraproduct over
$\mathcal{U}$ of the Hilbert spaces $H_{n}$. Following Raynaud (see
\cite{Ray}), we define $\mathcal{A}$, the vN-ultraproduct over
$\mathcal{U}$ of the von Neumann algebras $\mathcal{A}_{n}$, as the
$w^{*}$-closure of $\widetilde{\mathcal{A}}$ in $B(H)$. Then the
predual $\mathcal{A}_{*}$ of $\mathcal{A}$ is isometrically
isomorphic to the Banach ultraproduct over $\mathcal{U}$ of the
preduals $\left(\mathcal{A}_{n}\right)_{*}$:
\begin{equation}\label{pred}
\mathcal{A}_{*}=\Prod_{n\in
\N}\left(\mathcal{A}_{n}\right)_{*}/\mathcal{U}
\end{equation}
 Let us denote by $\varphi$ the normal state on $\mathcal{A}$ associated
to $(\varphi_{n})_{n\in \N}$. Note that $\varphi$ is not faithful on
$\mathcal{A}$, so we introduce $p \in \mathcal{A}$ the support of
the state $\varphi$. Recall that for all $x\in \mathcal{A}$ we have
$\varphi(x)=\varphi(xp)=\varphi(px)$, and that $\varphi(x)=0$ for a
positive $x$ implies that $pxp=0$. Denote by $(p\mathcal{A}p,
\varphi)$ the induced von Neumann algebra $p\mathcal{A}p\subset
B(pH)$ equipped with the restriction of the state $\varphi$. For
each $n \in \N$, let $(\sigma_{t}^{n})_{t\in \R}$ be the modular
group of automorphisms of $\varphi_{n}$ with the associated modular
operator given by $\Delta_{n}$. For all $t \in \R$, let
$(\Delta_{n}^{it})^{\bullet}$ be the associated unitary in $\Prod_{n
\in \N}B(H_{n})/\mathcal{U}\subset B(H)$. Since
$(\sigma_{t}^{n})_{n\in \N}^{\bullet}$ is the conjugation by
$(\Delta_{n}^{it})^{\bullet}$, it follows that
$(\sigma_{t}^{n})_{n\in \N}^{\bullet}$ extends by $w^{*}$-continuity
to a group of $*$-automorphisms of $\mathcal{A}$. Let
$(\sigma_{t})_{t \in \R}$ be the local modular group of
automorphisms of $p\mathcal{A}p$. By Raynaud's result (see Theorem
2.1 in \cite{Ray}), $p\mathcal{A}p$ is stable by
$(\sigma_{t}^{n})_{n\in \N}^{\bullet}$ and the restriction of
$(\sigma_{t}^{n})_{n\in \N}^{\bullet}$ to $p\mathcal{A}p$ coincides
with $\sigma_{t}$.

In the following, we consider  a von Neumann algebra
$\mathcal{N}\subset B(K)$ equipped with a normal faithful state
$\psi$. Let $\widetilde{\mathcal{N}}$ be a $w^{*}$-dense
$*$-subalgebra of $\mathcal{N}$ and $\Phi$ a $*$-homomorphism from
$\widetilde{\mathcal{N}}$ into $\mathcal{A}$ whose image will be
denoted by $\widetilde{\mathcal{B}}$ with $w^*$-closure denoted by
$\mathcal{B}$:
$$\Phi\;:\;\widetilde{\mathcal{N}}\subset\mathcal{N}\subset B(K)
\longrightarrow \widetilde{\mathcal{B}}\subset \mathcal{A} \subset
B(H)\quad\text{and}\quad
\overline{\widetilde{\mathcal{N}}}^{w^{*}}=\mathcal{N},\quad
\overline{\widetilde{\mathcal{B}}}^{w^{*}}=\mathcal{B}$$

By a result of Takesaki (see \cite{Tak}) there is a normal
conditional expectation from $p\mathcal{A}p$ onto $p\mathcal{B}p$ if
and only if $p\mathcal{B}p$ is stable by the modular group of
$\varphi$ (which is here given by Raynaud's results). Under this
condition there will be a normal conditional expectation from
$\mathcal{A}$ onto $p\mathcal{B}p$ and $p\mathcal{B}p$ will inherit
some of the properties of $\mathcal{A}$. We would like to pull back
these properties to $\mathcal{N}$ itself. It turns out that, with
good assumptions on $\Phi$ (see Lemma \ref{p} below), the
compression from $\mathcal{B}$ onto $p\mathcal{B}p$ is a
$*$-homomorphism. If in addition, we suppose that $\Phi$ is state
preserving, then $p\Phi p$ can be extended into a $w^{*}$-continuous
$*$-isomorphism between $\mathcal{N}$ and $p\mathcal{B}p$.

\begin{lem}\label{p} In the following,
$1.\Longrightarrow 2. \Longrightarrow 3. \Longleftrightarrow 4.
\Longleftrightarrow 5.$:
\begin{enumerate}
\item For all $x\in \widetilde{\mathcal{B}}$ there is a representative
$(x_{n})_{n \in \N}$ of $x$ such that for all $ n\in \N$, $x_{n}$ is
entire for $(\sigma_{t}^n)_{t \in \R}$ and
$(\sigma_{-i}^{n}(x_{n}))_{n \in \N}$ is uniformly bounded.
\item For all $x\in \widetilde{\mathcal{B}}$ there exists $z\in
\mathcal{A}$ such that for all $y \in \mathcal{A}$ we have
$\varphi(xy)=\varphi(yz)$.
\item For all $(x,y)\in \mathcal{B}^{2}$:
$\varphi(xpy)=\varphi(xy)$
\item For all $(x,y)\in \mathcal{B}^{2}$,
$pxyp=pxpyp$, i.e the canonical application from $\mathcal{B}$ to
$p\mathcal{B}p$ is a $*$-homomorphism.
\item $p\in \mathcal{B}'$.
\end{enumerate}
\end{lem}

\proof $1.\Longrightarrow 2.$ Consider $x \in
\widetilde{\mathcal{B}}$ with a representative $(x_{n})_{n \in \N}$
such that for all $ n\in \N$, $x_{n}$ is entire for
$(\sigma_{t}^n)_{t \in \R}$ and $(\sigma_{-i}^{n}(x_{n}))_{n \in
\N}$ is uniformly bounded. Denote by $z\in \mathcal{A}$ the class
$(\sigma_{-i}^{n}(x_{n}))_{n \in \N}^{\bullet}$. By $w^{*}$-density
and continuity if suffices to consider an element $y$ in
$\widetilde{\mathcal{A}}$ with representative $(y_{n})_{n \in \N}$.
Then,
\begin{equation*}
\varphi(xy) =\lim_{n,\mathcal{U}}\varphi_{n}(x_{n}
y_{n})=\lim_{n,\mathcal{U}}\varphi_{n}(y_{n}\sigma_{-i}^{n}(x_{n}))
=\varphi(yz)
\end{equation*}
$2.\Longrightarrow 3.$ Here again it suffices to consider $(x,y) \in
\widetilde{\mathcal{B}}^{2}$. By assumption there exists $z\in
\mathcal{A}$ such that for all $t\in \mathcal{A}$,
$\varphi(xt)=\varphi(tz)$. Applying our assumption for $t=py$ and
$t=y$ successively, we obtain the desired result:
\begin{equation*}
\varphi(xpy)=\varphi(pyz)=\varphi(yz)=\varphi(xy)
\end{equation*}
$3.\Longrightarrow 4.$ Let $x\in \mathcal{B}$. We have, by 3.:
$\varphi(x(1-p)x^{*})=0$. Since $p$ is the support of $\varphi$ and
$x(1-p)x^{*}\ge 0$, this implies $px(1-p)x^{*}p=0$. Thus for all
$x\in \mathcal{B}$ we have
$$pxpx^{*}p=pxx^{*}p$$
We conclude by polarization.\\
$4.\Longrightarrow 5.$ Let $q$ be an orthogonal projection in
$\mathcal{B}$. By $4.$, $pqp$ is again an orthogonal projection and
we claim that this is equivalent to $pq=qp$.
Indeed, let us denote by $x$ the contraction $qp$. Then $x^{*}x=pqp$
and since $pqp$ is an orthogonal projection we have $|x|=pqp$. It
follows that the polar decomposition of $x$ is of the form $x=upqp$,
with $u$ a partial isometry. Computing $x^{2}$, we see that $x$ is a
projection:
$$x^{2}=upqp(qp)=upqp=x.$$
Since $x$ is contractive, we deduce that $x$ is an orthogonal
projection and that $x^{*}=x$. Thus $pq=qp$. Since $\mathcal{B}$ is
generated by its
projections, we have $p\in \mathcal{B}'$.\\
$5.\Longrightarrow 3.$ This is clear.\Eproof

\medskip

We assume that one of the technical conditions of the previous Lemma
is fulfilled. Let us denote by $\Theta=p\Phi p$. $\Theta$ is a
$*$-homomorphism from $\widetilde{\mathcal{N}}$, into
$p\mathcal{A}p$.
$$\Theta=p\Phi p\;:\;\widetilde{\mathcal{N}}\longrightarrow
p\mathcal{A}p\subset B(pH)$$
 We assume that $\Phi$, and hence $\Theta$, is
state preserving. Then $\Theta$ can be extended into a
($w^{*}$-continuous) $*$-isomorphism from $\mathcal{N}$ onto
$p\mathcal{B}p$. This is indeed a consequence of the following well
known fact:

\begin{lem}\label{ext} Let $(\mathcal{M},\varphi)$ and $(\mathcal{N},\psi)$ be
von Neumann algebras equipped with normal faithful states. Let
$\widetilde{\mathcal{M}}$, (respectively $\widetilde{\mathcal{N}}$),
be a $w^{*}$ dense $*$-subalgebra of $\mathcal{M}$ (respectively
$\mathcal{N}$). Let $\Psi$ be a $*$-homomorphism from
$\widetilde{\mathcal{M}}$ onto $\widetilde{\mathcal{N}}$ such that
for all $m\in \widetilde{\mathcal{M}}$ we have
$\psi(\Psi(m))=\varphi(m)$ ($\Psi$ is state preserving). Then $\Psi$
extends uniquely into a normal $*$-isomorphism between $\mathcal{M}$
and $\mathcal{N}$.
\end{lem}

\proof Since $\varphi$ is faithful, we have for all $m\in
\mathcal{M}$, $\|m\|=\lim\limits_{n \to
+\infty}\varphi\left((m^{*}m)^{n}\right)^{\f{1}{2n}}$. Thus, since
$\Psi$ is state preserving, $\Psi$ is isometric from
$\widetilde{\mathcal{M}}$ onto $\widetilde{\mathcal{N}}$. We put
$$\varphi \widetilde{\mathcal{M}}=\{\varphi.m,\;m\in \widetilde{\mathcal{M}}\}
\subset \mathcal{M}_{*}\quad\text{and}\quad\psi
\widetilde{\mathcal{N}}=\{\psi.n,\;n\in
\widetilde{\mathcal{N}}\}\subset \mathcal{N}_{*}.$$
 $\varphi
\widetilde{\mathcal{M}}$ (respectively $\psi
\widetilde{\mathcal{N}}$) is dense in $\mathcal{M}_{*}$
(respectively $\mathcal{N}_{*}$). Let us define the following linear
operator $\Xi$ from $\psi \widetilde{\mathcal{N}}$ onto $\varphi
\widetilde{\mathcal{M}}$:
$$\Xi(\psi.\Psi(m))=\varphi.m\quad\text{for all\;} m\in
\widetilde{\mathcal{M}}$$ Using Kaplansky's density Theorem and the
fact that $\Psi$ is isometric, we compute:
\begin{align*}
\|\Xi(\psi .\Psi(m)\|&=\sup_{m_{0}\in
\widetilde{\mathcal{M}},\;\|m_{0}\|\le
1}\|\varphi(mm_{0})\|=\sup_{m_{0}\in
\widetilde{\mathcal{M}},\;\|m_{0}\|\le 1}\|\psi(\Psi(m)\Psi(m_{0}))\|\\
&=\sup_{n_{0}\in \widetilde{\mathcal{N}},\;\|n_{0}\|\le
1}\|\psi(\Psi(m)n_{0}))\|=\|\psi.\Psi(m)\|
\end{align*}
So that $\Xi$ extends into a surjective isometry from
$\mathcal{N}_{*}$ onto $\mathcal{M}_{*}$. Moreover $\Xi$ is the
preadjoint of $\Psi$. Indeed we have for all $(m,m_{0})\in
\widetilde{\mathcal{M}}^{2}$:
$$\scal{\psi.\Psi(m)}{\Psi(m_{0})}=\psi(\Psi(m)\Psi(m_{0}))=
\varphi(mm_{0})=\scal{\Xi(\psi.\Psi(m))}{m_{0}}$$ Thus $\Psi$
extends to a normal $*$-isomorphism between $\mathcal{N}$ and
$\mathcal{M}$. \Eproof

\medskip

In the following Theorem, we sum up what we have proved in the
previous discussion:
\begin{thm}\label{qwep} Let $(\mathcal{N},\psi)$ and $(\mathcal{A}_{n},
\varphi_{n})$, for $n \in \N$, be von Neumann algebras equipped with
normal faithful states. Let $\mathcal{U}$ be a non trivial
ultrafilter on $\N$, and $\mathcal{A}$ the von Neumann algebra
ultraproduct over $\mathcal{U}$ of the $\mathcal{A}_{n}$'s. For all
$n \in \N$ let us denote by $(\sigma^{n}_{t})_{t\in \R}$ the modular
group of $\varphi_{n}$ and by $\varphi$ the normal state on
$\mathcal{A}$ which is the ultraproduct of the states $\varphi_{n}$.
$p\in \mathcal{A}$ denote the support of $\varphi$. Consider
$\widetilde{\mathcal{N}}$ a $w^{*}$-dense $*$-subalgebra of
$\mathcal{N}$ and  a $*$-homomorphism $\Phi$
$$\Phi\;:\;\widetilde{\mathcal{N}}\subset \mathcal{N}\longrightarrow
\mathcal{A}=\Prod_{n, \mathcal{U}}\mathcal{A}_{n}$$ Assume $\Phi$
satisfies:
\begin{enumerate}
\item $\Phi$ is state preserving: for all $x\in
\widetilde{\mathcal{N}}$ we have
$$\varphi(\Phi(x))=\psi(x)$$
\item For all $(x,y) \in \Phi(\widetilde{\mathcal{N}})^{2}$
$$\varphi(xy)=\varphi(xpy).$$
(Or one of the technical conditions of Lemma \ref{p}.)
\item For all $t\in \R$ and for all $y=(y_{n})_{n\in \N}^{\bullet}\in \Phi(\widetilde{\mathcal{N}})$,
$$p(\sigma_{t}^{n}(y_{n}))_{n\in \N}^{\bullet}p\;(=\sigma_{t}(pyp))\in p\mathcal{B}p$$
where $\mathcal{B}$ is the $w^{*}$-closure of
$\Phi(\widetilde{\mathcal{N}})$ in $\mathcal{A}$.
\end{enumerate}
Then $\Theta=p\Phi p\;:\;\widetilde{\mathcal{N}}\longrightarrow
p\mathcal{A}p$ is a state preserving $*$-homomorphism which can be
extended into a normal isomorphism (still denoted by $\Theta$)
between $\mathcal{N}$ and its image
$\Theta(\mathcal{N})=p\mathcal{B}p$. Moreover there exists a
(normal) state preserving conditional expectation from $\mathcal{A}$
onto $\Theta(\mathcal{N})$.
\end{thm}

\bigskip

{\noindent \it Remarks.}

\smallskip

{\noindent $\bullet$} Condition $2.$ is in fact necessary for
$\Theta$ being a $*$-homomorphism (by Lemma \ref{p}), and condition
$3.$ is necessary for the existence of a state preserving
conditional expectation onto $\Theta(\mathcal{N})$ (by \cite{Tak}).

\smallskip

{\noindent $\bullet$} Let us denote by $(\sigma_{t}^{\psi})_{t \in
\R}$ the modular group of $*$-automorphisms of $\psi$. Provided that
$(\sigma_{t}^{\psi})_{t \in \R}$ maps $\widetilde{\mathcal{N}}$ into
itself, we can replace condition 2. of the previous Theorem by the
following intertwining condition: For all $t\in \R$ and for all
$x\in \widetilde{\mathcal{N}}$ we have
$$p(\sigma_{t}^{n}(y_{n}))_{n\in \N}^{\bullet}p\;(=
\sigma_{t}(p\Phi(x)p))= p\Phi(\sigma_{t}^{\psi}(x))p$$ where
$\Phi(x)=(y_{n})_{n\in \N}^{\bullet}$. Moreover, notice that if the
conclusion of the Theorem is true, then this condition must be
fulfilled for all $t\in \R$ and for all $x\in \mathcal{N}$ (see
\cite{Tak2} page 95).

\begin{cor}\label{corqwep} Under the assumptions of the previous Theorem,
$\mathcal{N}$ is QWEP provided that each of the $\mathcal{A}_n$ is
QWEP.
\end{cor}

\proof This is a consequence of Kirchberg's results (see
\cite{Kir,Oz}). First, $\Prod_{n\in \N}\mathcal{A}_n$ is QWEP as a
product of QWEP $C^{*}$-algebras (\cite{Oz} Proposition 4.1 (i)).
Since $\tilde{\mathcal{A}}$ is a quotient of a QWEP $C^{*}$-algebra,
it is also QWEP. It follows that $\mathcal{A}$ which is the
$w^{*}$-closure of $\tilde{\mathcal{A}}$ in $B(H)$ is QWEP (by
\cite{Oz} Proposition 4.1 (iii)). Since there is a conditional
expectation from $\mathcal{A}$ onto $p\mathcal{A}p$, $p\mathcal{A}p$
is QWEP (see \cite{Kir}). Finally, by Theorem \ref{qwep},
$\mathcal{N}$ is isomorphic to a subalgebra of $p\mathcal{A}p$ which
is the image of a (state preserving) conditional expectation, thus
$\mathcal{N}$ inherits the QWEP property. \Eproof

\section{The finite dimensional case}\label{section4}

In this section we show that $\Gamma_{q}(H_{\R},(U_{t})_{t \in \R})$
is QWEP when $H_{\R}$ is finite dimensional. For notational purpose,
it will be more convenient to deal with $\dim (H_{\R})$ even. This
is not relevant in our context (see the remark after Theorem
\ref{qwepf}). We put $\dim (H_{\R})=2k$. Notice that
$\Gamma_{q}(H_{\R},(U_{t})_{t \in \R})$ only depends on the spectrum
of the operator $A$. The spectrum of $A$ is given by the set
$\{\lambda_{1},\dots, \lambda_{k}\}\cup \{\lambda_{1}^{-1},\dots,
\lambda_{k}^{-1}\}$ where for all $j \in \{1,\dots,k\}$,
$\lambda_{j}\ge 1$. As in subsection \ref{subsec2} , we use the
notation $\mu_{j}=\lambda_{j}^{\f{1}{4}}$.

\subsection{Twisted Baby Fock}

We start by adapting Biane's model to our situation. Let us denote
by $I$ the set $\{-k,\dots,-1\}\cup \{1,\dots,k\}$. As in subsection
\ref{subsec4}, we give us a function $\epsilon$ on $I\times I$ into
$\{-1,1\}$ and we consider the associated complex $*$-algebra
$\mathcal{A}(I,\epsilon)$. By analogy with (\ref{semicirc}), for all
$j \in \{1, \dots, k\}$ we define the following generalized
semi-circular variables acting on $L^{2}(\mathcal{A}(I,\epsilon),
\varphi^{\epsilon})$:
$$\gamma_{i}=\mu_{i}^{-1}\beta_{i}^{*}+\mu_{i}\beta_{-i}\quad{\text and}\quad\delta_{i}=\mu_{i}
\alpha_{i}^{*}+\mu_{i}^{-1}\alpha_{-i}$$ We denote by $\Gamma$
(respectively $\Gamma_{r}$) the von Neumann algebra generated in
$B(L^{2}(\mathcal{A}(I,\epsilon),\varphi^{\epsilon}))$ by the
$\gamma_{i}$ (respectively $\delta_{i}$). $\Gamma_{r}$ is the
natural candidate for the commutant of $\Gamma$ in
\\$B(L^{2}(\mathcal{A}(I,\epsilon),\varphi^{\epsilon}))$. We need to show that the vector
$1$ is cyclic and separating for $\Gamma$. To do so we must assume
that $\epsilon$ satisfies the following additional condition:
\begin{equation}\label{abso}
\text{For all\;}\quad (i,j) \in I^{2},\quad
\epsilon(i,j)=\epsilon(|i|,|j|)
\end{equation}
This condition is in fact a necessary condition for
$\Gamma_{r}\subset \Gamma '$ and for condition 1.(a) of Lemma
\ref{rep} below.
\begin{lem}\label{plus} Under condition (\ref{abso}) the following relation
holds: $$\text{For all}\quad i\in I,\quad
\alpha_{i}\beta_{i}^{*}+\alpha_{-i}^{*}\beta_{-i}=\beta_{i}^{*}\alpha_{i}
+\beta_{-i}\alpha_{-i}^{*}$$
\end{lem}

\proof Let $i \in I$ and $A\subset I$. We have
$$(\alpha_{i}\beta_{i}^{*}+\alpha_{-i}^{*}\beta_{-i})(x_{A})=
\left\{ \begin{array}{cl}
x_{-i}x_{A}x_{-i} & {\rm if}\quad i\in A\quad{\rm and}\quad -i\in A\\
0 & {\rm if}\quad i\in A\quad{\rm and}\quad -i\not\in A\\
x_{i}x_{A}x_{i}+x_{-i}x_{A}x_{-i} & {\rm if}\quad i\not\in A\quad{\rm and}\quad -i\in A\\
x_{i}x_{A}x_{i}& {\rm if}\quad i\not\in A\quad{\rm and}\quad -i\not\in A\\
\end{array}
\right.$$ and
$$(\beta_{i}^{*}\alpha_{i}+\beta_{-i}\alpha_{-i}^{*})(x_{A})=
\left\{ \begin{array}{cl}
x_{i}x_{A}x_{i} & {\rm if}\quad i\in A\quad{\rm and}\quad -i\in A\\
x_{i}x_{A}x_{i}+x_{-i}x_{A}x_{-i} & {\rm if}\quad i\in A\quad{\rm and}\quad -i\not\in A\\
0 & {\rm if}\quad i\not\in A\quad{\rm and}\quad -i\in A\\
x_{-i}x_{A}x_{-i}& {\rm if}\quad i\not\in A\quad{\rm and}\quad -i\not\in A\\
\end{array}
\right.$$

Thus, we need to study the following cases. Assume that
$A=\{i_{1},\dots,i_{p}\}$ where $i_{1}<\dots<i_{p}$.
\begin{enumerate}
\item If $i$ and $-i$ belong to $A$ then there exists $(l,m) \in \{1,\dots,p\}$, $l<m$,  such that $i_{l}=-i$
and $i_{m}=i$. Applying successively relations (\ref{mix}) and
(\ref{abso}), we get:
\begin{align*}
x_{-i}x_{A}x_{-i} = & \left(\prod_{q=1}^{l-1}\epsilon(i_{q},-i)
\right)
x_{i_{1}}\dots x_{i_{l-1}}x_{i_{l+1}}\dots x_{i_{p}}x_{-i}\\
= & \left(\prod_{q=1}^{l-1}\epsilon(i_{q},-i)\right)
\left(\prod_{q=l+1}^{p}\epsilon(i_{q},-i)\right)x_{A}=-\left(\prod_{q=1}^{p}\epsilon(i_{q},-i)\right)x_{A}\\
= & -\left(\prod_{q=1}^{p}\epsilon(i_{q},i)\right)x_{A}=
x_{i}x_{A}x_{i}
\end{align*}
\item If $i$ and $-i$ do not belong to $A$, we can check in a similar way that:
\begin{align*}
x_{-i}x_{A}x_{-i} = &
\left(\prod_{q=1}^{p}\epsilon(i_{q},-i)\right)x_{A}=
\left(\prod_{q=1}^{p}\epsilon(i_{q},i)\right)
x_{A}\\
= & x_{i}x_{A}x_{i}
\end{align*}
\item If $i\in A$ and $-i \not\in A$, then there exists $l\in \{1,\dots,p\}$ such that $i_{l}=i$. We have:
\begin{align*}
x_{i}x_{A}x_{i} = & \left(\prod_{q=1}^{l-1}\epsilon(i_{q},i)
\right)x_{i_{1}}\dots x_{i_{l-1}}
x_{i_{l+1}}\dots x_{i_{p}}x_{i}\\
= & \left(\prod_{q=1}^{l-1}\epsilon(i_{q},i)
\right)\left(\prod_{q=l+1}^{p}\epsilon(i_{q},i)\right)x_{A}
=-\left(\prod_{q=1}^{p}\epsilon(i_{q},i)\right)x_{A}\\
= &
-\left(\prod_{q=1}^{p}\epsilon(i_{q},-i)\right)x_{A}=-x_{-i}x_{A}x_{-i}
\end{align*}
This finishes the proof.
\end{enumerate}
\Eproof

\begin{lem}\label{rep} By construction we have:
\begin{enumerate}
\item For all $(i,j)\in \{1, \dots, k\}^{2}$, $i \neq j$, the
following mixed commutation and anti-commutation relations hold:
\begin{enumerate}
\item $\gamma_{i}\gamma_{j}-\epsilon(i,j)\gamma_{j}\gamma_{i}=0$
\item $\gamma^{*}_{i}\gamma_{j}-\epsilon(i,j)\gamma_{j}\gamma^{*}_{i}=0$
\item $(\gamma_{i}^{*})^{2}=\gamma_{i}^{2}=0$
\item $\gamma_{i}^{*}\gamma_{i}+\gamma_{i}\gamma_{i}^{*}=(\mu_{i}^{2}+\mu_{i}^{-2})Id$.
\end{enumerate}
\item Same relations as in 1.  for the operators $\delta_{i}$.
\item $\Gamma_{r}\subset \Gamma'$. \item The vector $1$ is cyclic
and separating for both $\Gamma$ and $\Gamma_{r}$. \item $\Gamma
\subset B(L^{2}(\mathcal{A}(I,\epsilon),\varphi^{\epsilon}))$ is the
(faithful) G.N.S representation of
$\left(\Gamma,\varphi^{\epsilon}\right)$.
\end{enumerate}
\end{lem}

\proof 1.(a) Thanks to 2. of Lemma \ref{relcrea} and (\ref{abso}) we
get:
\begin{align*}
\gamma_{i}\gamma_{j}  &=
\mu_{i}^{-1}\mu_{j}^{-1}\beta_{i}^{*}\beta_{j}^{*}+
\mu_{i}\mu_{j}\beta_{-i}\beta_{-j}+\mu_{i}^{-1}\mu_{j}\beta_{i}^{*}\beta_{-j}+
\mu_{i}\mu_{j}^{-1}\beta_{-i}\beta_{j}^{*}\\
&= \epsilon(i,j)\mu_{i}^{-1}\mu_{j}^{-1}\beta_{j}^{*}\beta_{i}^{*}
+\epsilon(-i,-j)\mu_{i}\mu_{j}\beta_{-j}\beta_{-i}+
\epsilon(i,-j)\mu_{i}^{-1}\mu_{j}\beta_{-j}\beta_{i}^{*}\\
&
\phantom{=\epsilon(i,j)\mu_{i}^{-1}\mu_{j}^{-1}\beta_{j}^{*}\beta_{i}^{*}\,\,}
+\epsilon(-i,j)\mu_{i}\mu_{j}^{-1}\beta_{j}^{*}\beta_{-i}\\
&=
\epsilon(i,j)\left(\mu_{i}^{-1}\mu_{j}^{-1}\beta_{j}^{*}\beta_{i}^{*}
+
\mu_{i}\mu_{j}\beta_{-j}\beta_{-i}+\mu_{i}^{-1}\mu_{j}\beta_{-j}\beta_{i}^{*}
+\mu_{i}\mu_{j}^{-1}\beta_{j}^{*}\beta_{-i}\right)\\
&= \epsilon(i,j)\gamma_{j}\gamma_{i}
\end{align*}
1.(b) Is analogous to (a) and is left to the reader.\\
1.(c) Using 1. and 2. of Lemma \ref{relcrea}, and
$\epsilon(i,-i)=\epsilon(i,i)=-1$ we get:
\begin{align*}
\gamma_{i}^{2}&
=\mu_{i}^{-2}(\beta_{i}^{*})^{2}+\mu_{i}^{2}\beta_{-i}^{2}+\beta_{i}^{*}\beta_{-i}+\beta_{-i}
\beta_{i}^{*}\\
& =\epsilon(i,-i)\beta_{-i}\beta_{i}^{*}+\beta_{-i}\beta_{i}^{*}=0
\end{align*}
1.(d) Using similar arguments, we compute:
\begin{align*}
\gamma_{i}^{*}\gamma_{i}+\gamma_{i}\gamma_{i}^{*}
&=\mu_{i}^{-2}(\beta_{i}\beta_{i}^{*}+\beta_{i}^{*}\beta_{i})
+\mu_{i}^{2}(\beta_{-i}^{*}\beta_{-i}+\beta_{-i}\beta_{-i}^{*})+\beta_{i}\beta_{-i}+\beta_{-i}\beta_{i}\\
&\phantom{=\mu_{i}^{-2}(\beta_{i}\beta_{i}^{*}+\beta_{i}^{*}\beta_{i})\,\,}
+\beta_{-i}^{*}\beta_{i}^{*}+\beta_{i}^{*}\beta_{-i}^{*}\\
&=(\mu_{i}^{-2}+\mu_{i}^{2})Id+(\epsilon(i,-i)+1)(\beta_{i}\beta_{-i}+\beta_{-i}^{*}\beta_{i}^{*})
=(\mu_{i}^{-2}+\mu_{i}^{2})Id
\end{align*}
2. Is now clear from the proof of 1. since the relations for the
$\alpha_{i}$'s are the same as the ones for the $\beta_{i}$'s.\\
3. It suffices to show that for all $(i,j) \in \{1,\dots,k\}^{2}$ we
have $\gamma_{i}\delta_{j}=\delta_{j}\gamma_{i}$ and
$\gamma_{i}\delta_{j}^{*}=\delta_{j}^{*}\gamma_{i}$.\\
If $i\neq j$ then from 5. of Lemma \ref{relcrea} it is clear that
$\gamma_{i}\delta_{j}=\delta_{j}\gamma_{i}$ and
$\gamma_{i}\delta_{j}^{*}=\delta_{j}^{*}\gamma_{i}$.\\
If $i=j$ then using 4. and 5. of Lemma \ref{relcrea} and Lemma
\ref{plus} we obtain the desired result as follows:
\begin{align*}
\gamma_{i}\delta_{i} &=
\beta_{i}^{*}\alpha_{i}^{*}+\beta_{-i}\alpha_{-i}
+\mu_{i}^{-2}\beta_{i}^{*}\alpha_{-i}+\mu_{i}^{2}\beta_{-i}\alpha_{i}^{*}
=\mu_{i}^{-2}\beta_{i}^{*}\alpha_{-i}+\mu_{i}^{2}\beta_{-i}\alpha_{i}^{*}\\
&=\mu_{i}^{-2}\alpha_{-i}\beta_{i}^{*}+\mu_{i}^{2}\alpha_{i}^{*}\beta_{-i}=\delta_{i}\gamma_{i}
\end{align*}
and
\begin{align*}
\gamma_{i}\delta_{i}^{*}
&=\beta_{i}^{*}\alpha_{i}+\beta_{-i}\alpha_{-i}^{*}
+\mu_{i}^{-2}\beta_{i}^{*}\alpha_{-i}^{*}+\mu_{i}^{2}\beta_{-i}\alpha_{i}\\
&=\alpha_{i}\beta_{i}^{*}+\alpha_{-i}^{*}\beta_{-i}+\mu_{i}^{-2}\beta_{i}^{*}\alpha_{-i}
+\mu_{i}^{2}\beta_{-i}\alpha_{i}^{*}\\
&=\alpha_{i}\beta_{i}^{*}+\alpha_{-i}^{*}\beta_{-i}+\mu_{i}^{-2}\alpha_{-i}\beta_{i}^{*}
+\mu_{i}^{2}\alpha_{i}^{*}\beta_{-i}=\delta_{i}^{*}\gamma_{i}
\end{align*}
4. It suffices to prove that for any $A \subset I$ we have $x_{A}\in
\Gamma 1\cap \Gamma_{r}1$. Let $A\subset I$ and $(\chi_{i})_{i\in
I}\in \{0,1\}^{I}$ such that $\chi_{i}=1$ if and only if $i \in A$.
Then
\begin{align*}
x_{A} &=x_{-k}^{\chi_{-k}}\dots
x_{-1}^{\chi_{-1}}x_{1}^{\chi_{1}}\dots x_{k}^{\chi_{k}}\\
&= (\mu_{k}^{-1}\gamma_{k}^{*})^{\chi_{-k}}\dots
(\mu_{1}^{-1}\gamma_{1}^{*})^{\chi_{-1}}(\mu_{1}\gamma_{1})^{\chi_{1}}\dots
(\mu_{k}\gamma_{k})^{\chi_{k}}1\\
&= \mu_{1}^{\chi_{1}-\chi_{-1}}\dots
\mu_{k}^{\chi_{k}-\chi_{-k}}\gamma_{k}^{-\chi_{-k}}\dots
\gamma_{1}^{-\chi_{-1}}\gamma_{1}^{\chi_{1}}\dots
\gamma_{k}^{\chi_{k}}1
\end{align*}
where by convention $\gamma_{i}^{-1}=\gamma_{i}^{*}$. \\
The same computation is valid for $\Gamma_{r}$ and we obtain:
$$x_{A}=\mu_{1}^{\chi_{-1}-\chi_{1}}\dots
\mu_{k}^{\chi_{-k}-\chi_{k}}\delta_{k}^{\chi_{k}}\dots
\delta_{1}^{\chi_{1}}\delta_{1}^{-\chi_{-1}}\dots
\delta_{k}^{-\chi_{-k}}1$$ It follows that the vector $1$ is cyclic
for both $\Gamma$ and $\Gamma_{r}$. Since $\Gamma_{r}\subset \Gamma
'$ then $1$ is also cyclic for $\Gamma '$ and thus separating for
$\Gamma$. The same argument applies to $\Gamma_{r}$ and thus $1$ is
also a cyclic and separating vector for $\Gamma_{r}$.\\
5. This is clear from the just proved assertion and the fact that
the state $\varphi^{\epsilon}$ is equal to the vector state
associated to the vector $1$. \Eproof

By the Lemma just proved, we are in a situation where we can apply
Tomita-Takesaki theory. As usual we denote by $S$ the involution on
$L^{2}(\mathcal{A}(I,\epsilon),\varphi^{\epsilon})$ defined by:
$S(\gamma 1)=\gamma^{*}1$ for all $\gamma \in \Gamma$. $\Delta$ will
denote the modular operator and $J$ the modular conjugation. Recall
that $S=J\Delta^{\f{1}{2}}$ is the polar decomposition of the
antilinear operator $S$ (which is here bounded since we are in a
finite dimensional framework). We also denote by $(\sigma_{t})_{t\in
\R}$ the modular group of automorphisms of $\Gamma$ associated to
$\varphi$. Recall that for all $\gamma \in \Gamma$ and all $t \in
\R$ we have $\sigma_{t}(\gamma)=\Delta^{it}\gamma\Delta^{-it}$.

\medskip

{\noindent Notation}: In the following, for $A\subset I$ we denote
by $(\chi_{i})_{i \in I}$ the characteristic function of the set $A$
: $\chi_{i}=1$ if $i \in A$ and $\chi_{i}=0$ if $i \not\in A$. (We
will not keep track of the dependance in $A$ unless there could be
some confusion.)

\begin{prop}\label{modu}The modular operators and the modular group of $(\Gamma, \varphi^{\epsilon})$
are determined by:
\begin{enumerate}
 \item $J$ is the antilinear operator given by: for all $A\subset
 I$,
 $$J(x_{A})=J(x_{-k}^{\chi_{-k}}\dots x_{-1}^{\chi_{-1}} x_{1}^{\chi_{1}}\dots x_{k}^{\chi_{k}})=
 x_{-k}^{\chi_{k}}\dots x_{-1}^{\chi_{1}} x_{1}^{\chi_{-1}}\dots
 x_{k}^{\chi_{-k}}$$
 \item $\Delta$ is the diagonal and positive operator given by: for all
 $A\subset I$,
 $$\Delta(x_{A})=\Delta(x_{-k}^{\chi_{-k}}\dots x_{-1}^{\chi_{-1}} x_{1}^{\chi_{1}}\dots
 x_{k}^{\chi_{k}})=\lambda_{k}^{(\chi_{k}-\chi_{-k})}\dots
 \lambda_{1}^{(\chi_{1}-\chi_{-1})}x_{A}$$
 \item For all $j \in \{1\dots ,k\}$, $\gamma_{j}$ is entire for $(\sigma_{t})_{t}$ and
 satisfies $\sigma_{z}(\gamma_{j})=\lambda_{j}^{iz}\gamma_{j}$ for
 all $z \in \C$.
\end{enumerate}
\end{prop}

\proof Let $A\subset I$. We have
$$x_{A}=x_{-k}^{\chi_{-k}}\dots
x_{-1}^{\chi_{-1}} x_{1}^{\chi_{1}}\dots
x_{k}^{\chi_{k}}=\mu_{1}^{\chi_{1}-\chi_{-1}}\dots
\mu_{k}^{\chi_{k}-\chi_{-k}}\gamma_{k}^{-\chi_{-k}}\dots
\gamma_{1}^{-\chi_{-1}}\gamma_{1}^{\chi_{1}}\dots
\gamma_{k}^{\chi_{k}}1$$ Thus,
\begin{align*}
S(x_{A}) &=\mu_{1}^{\chi_{1}-\chi_{-1}}\dots
\mu_{k}^{\chi_{k}-\chi_{-k}}(\gamma_{k}^{-\chi_{-k}}\dots
\gamma_{1}^{-\chi_{-1}}\gamma_{1}^{\chi_{1}}\dots
\gamma_{k}^{\chi_{k}})^{*}1\\
&=\mu_{1}^{\chi_{1}-\chi_{-1}}\dots
\mu_{k}^{\chi_{k}-\chi_{-k}}\gamma_{k}^{-\chi_{k}}\dots
\gamma_{1}^{-\chi_{1}}\gamma_{1}^{\chi_{-1}}\dots\gamma_{k}^{\chi_{-k}}1\\
&=\mu_{1}^{2(\chi_{1}-\chi_{-1})}\dots
\mu_{k}^{2(\chi_{k}-\chi_{-k})}x_{-k}^{\chi_{k}}\dots
x_{-1}^{\chi_{1}} x_{1}^{\chi_{-1}}\dots
 x_{k}^{\chi_{-k}}
\end{align*}
By uniqueness of the polar decomposition, we obtain the stated
result. Let $j \in \{1\dots k\}$ and $t \in \R$ we have:
\begin{align*}
\sigma_{t}(\gamma_{j})1&
=\Delta^{it}\gamma_{j}\Delta^{-it}1=\Delta^{it}\gamma_{j}1
=\mu_{j}^{-1}\Delta^{it}x_{j}=\mu_{j}^{-1}\mu_{j}^{4it}x_{j}\\
&= \mu_{j}^{4it}\gamma_{j}1
\end{align*}
It follows, since $1$ is separating for $\Gamma$, that
$\sigma_{t}(\gamma_{j})=\mu_{j}^{4it}\gamma_{j}$. \Eproof

\medskip

{\noindent \it Remarks.}

\smallskip

{\noindent $\bullet$} We have $\Gamma '=\Gamma_{r}$. Indeed we have
already proved the inclusion $\Gamma_{r}\subset \Gamma '$ in Lemma
\ref{rep} . For the reverse inclusion we can use Tomita-Takesaki
theory which ensures that $\Gamma '= J\Gamma J$. But for all $j\in
I$ it is easy to see that $J\beta_{j}J=\alpha_{-j}$. It follows that
for all $j\in \{1,\dots,k\}$ we have $J\gamma_{j} J=\delta_{j}^{*}$.
Thus $\Gamma ' \subset \Gamma_{r}$. The equality $\Gamma
'=\Gamma_{r}$ can also be seen as a consequence of a general fact in
Tomita-Takesaki theory: it suffices to remark that $\Gamma_{r}$ is
the right Hilbertian algebra associated to $\Gamma$ in its GNS
representation.

\smallskip

{\noindent $\bullet$} The previous construction can be performed for
an infinite set of the form $J\times \{-1,1\}$ given with a family
of eigenvalues $(\mu_{j})_{j\in J} \in [1, +\infty[^{J}$ and a sign
function $\epsilon$ satisfying
$$\epsilon((j,i),(j^{\prime},i^{\prime}))=\epsilon((j,1),(j^{\prime},1))\quad\text{for all}\quad
((j,i),(j^{\prime},i^{\prime}))\in (J\times \{-1,1\})^{2}.$$

\subsection{Central limit approximation of $q$-Gaussians}

In this section we use the twisted Baby Fock construction to obtain
an asymptotic random matrix model for the $q$-Gaussian variables,
via Speicher's central limit Theorem. Let us first check the
independence condition:

\begin{lem}\label{ind} For all $j \in \{1,\dots,k\}$ let us denote by
$\mathcal{A}_{j}$ the $C^{*}$-subalgebra of\\
$B(L^{2}(\mathcal{A}(I,\epsilon),\varphi^{\epsilon}))$ generated by
the operators $\beta_{j}$ and $\beta_{-j}$. Then the family
$(\mathcal{A}_{j})_{1\le j \le k}$ is independent in
$B(L^{2}(\mathcal{A}(I,\epsilon),\varphi^{\epsilon}))$. In
particular, the family $(\gamma_{j})_{1\le j \le k}$ is independent.
\end{lem}

\proof The proof proceeds by induction. Changing notation, it
suffices to show that
$$\varphi^{\epsilon}(a_{1}\dots a_{r+1})=\varphi^{\epsilon}(a_{1}\dots
a_{r})\varphi^{\epsilon}(a_{r+1})$$ where $a_{l}\in \mathcal{A}_{l}$
for all $l\in \{1,\dots, r+1\}$. Since $a_{r+1}$ is a certain
non-commutative polynomial in the variables $\beta_{r+1}$,
$\beta^{*}_{r+1}$, $\beta_{-(r+1)}$, and $\beta^{*}_{-(r+1)}$, it is
clear that there exists $\nu \in \text{Span}\{x_{r+1},\;
x_{-(r+1)},\;x_{-(r+1)}x_{r+1}\}$ such that
$$a_{r+1}1=\scal{1}{a_{r+1}1}1+\nu$$
It is easy to see that $a_{r}^{*}\dots a_{1}^{*}1 \in
\text{Span}\left\{x_{B},\;B\subset
\{-r,\dots,-1\}\cup\{1,\dots,r\}\right\}$, which is orthogonal to
$\text{Span}\{x_{r+1},\; x_{-(r+1)},\;x_{-(r+1)}x_{r+1}\}$. We
compute:
\begin{align*}
\varphi^{\epsilon}(a_{1}\dots a_{r+1})& =\scal{1}{a_{1}\dots a_{r}
a_{r+1}1}=\scal{a_{r}^{*}\dots a_{1}^{*}1}{a_{r+1}1}\\
& =\scal{a_{r}^{*}\dots
a_{1}^{*}1}{1}\scal{1}{a_{r+1}1}+\scal{a_{r}^{*}\dots
a_{1}^{*}1}{\nu}=\scal{1}{a_{1}\dots
a_{r}1}\scal{1}{a_{r+1}1}\\
& = \varphi^{\epsilon}(a_{1}\dots a_{r})\varphi^{\epsilon}(a_{r+1})
\end{align*}
\Eproof

\bigskip

{\noindent \it Remark.} It is clear that one can prove, in the same
way, that the $C^*$-algebras generated by the $\beta_{j}$ are
independent (this is Proposition $3$ in \cite{Bi}).

\bigskip

Let $q \in (-1,1)$. Let us choose a family of random variables
$(\epsilon(i,j))_{(i,j)\in \N_{*}^{2},i\neq j}$ as in Lemma
\ref{rand}, and set $\epsilon(i,i)=-1$ for all $i \in \N_{*}$. As in
section 2.5.1, for all $n\in \N_{*}$ we will consider the complex
$*$-algebra $\mathcal{A}(I_{n},\epsilon_{n})$ where
$$I_{n}=\{1,\dots,n\}\times\big(\{-k,\dots,-1\}\cup\{1,\dots,k\}\big)$$ and
$$\epsilon_{n}((i,j),(i^{\prime},j^{\prime}))=\epsilon(i,i^{\prime})\quad \text{for all}
\quad ((i,j),(i^{\prime},j^{\prime}))\in I_{n}^{2}.$$ Notice that
the analogue of condition (\ref{abso}) is automatically satisfied.
Indeed, we have:
$$\epsilon_{n}((i,j),(i^{\prime},j^{\prime}))=\epsilon_{n}((i,|j|),(i^{\prime},|j^{\prime}|))
\quad \text{for all} \quad ((i,j),(i^{\prime},j^{\prime}))\in
I_{n}^{2}.$$
 Let us remind that $\mathcal{A}(I_{n},\epsilon_{n})$ is
the unital free complex algebra with generators $(x_{i,j})_{(i,j)\in
I_n}$ quotiented by the relations,
$$x_{i,j}x_{i^{\prime},j^{\prime}} -\epsilon(i,i^{\prime})x_{i^{\prime},j^{\prime}}x_{i,j}
=2\delta_{(i,j),(i^{\prime},j^{\prime})}$$ and with involution given
by $x_{i,j}^{*}=x_{i,j}$. For all $(i,j) \in
\{1,\dots,n\}\times\{1,\dots,k\}$ let $\gamma_{i,j}$ be the "twisted
semi-circular variable" associated to $\mu_{j}$
$$\gamma_{i,j}=\mu_{j}^{-1}\beta^{*}_{i,j}+\mu_{j}\beta_{i,j}$$
We denote by $\Gamma_{n}\subset B(L^{2}(\mathcal{A}(I_n,
\epsilon_n),\varphi^{\epsilon_n}))$ the von-Neumann algebra
generated by the $\gamma_{i,j}$ for $(i,j) \in
\{1,\dots,n\}\times\{1,\dots,k\}$. Observe that all our notations
are consistent since $(\Gamma_{n}, \varphi^{\epsilon_n})$ is
naturally embedded in $(\Gamma_{n+1}, \varphi^{\epsilon_{n+1}})$
(see the remarks following Lemma (\ref{relcrea})). In fact all these
algebras $(\Gamma_n, \varphi^{\epsilon_{n}})$ can be embedded in the
bigger von Neumann algebra $(\Gamma, \varphi^{\overline{\epsilon}})$
which is the Baby Fock construction associated to the infinite set
$\overline{I}$ and the sign function $\overline{\epsilon}$ given by
$$\overline{I}=\N_{*}\times \big(\{-k,\dots,-1\}\cup\{1,\dots,k\}\big) $$
and
$$\overline{\epsilon}((i,j),(i^{\prime},j^{\prime}))=\epsilon(i,i^{\prime})\quad \text{for all}
\quad ((i,j),(i^{\prime},j^{\prime}))\in \overline{I}^{2}.$$ Let us
denote by $s_{n,j}$ the following sum:
$$s_{n,j}=\f{1}{\sqrt{n}}\sum_{i=1}^{n}\gamma_{i,j}$$
We now check the hypothesis of Theorem \ref{TCL} for the family
$(\gamma_{i,j})_{(i,j)\in \N_{*}\times \{1,\dots,k\}}\subset
(\Gamma, \varphi^{\overline{\epsilon}})$.
\begin{enumerate}
\item The family is independent by Lemma \ref{ind}.
\item It is clear that for all $(i,j)$ we have
$\varphi^{\overline{\epsilon}}(\gamma_{i,j})=0$.
\item Let $(j(1),j(2))\in \{1,\dots,k\}$ and $i \in \N_{*}$. We
compute and identify the covariance thanks to Lemma \ref{momdef}:
\begin{align*}
\varphi^{\overline{\epsilon}}(\gamma_{i,j(1)}^{k(1)}\gamma_{i,j(2)}^{k(2)})&=\scal{\gamma_{i,j(1)}^{-k(1)}1}
{\gamma_{i,j(2)}^{k(2)}1}=
\scal{\mu_{j(1)}^{k(1)}x_{-k(1)i,-k(1)j(1)}}{\mu_{j(2)}^{-k(2)}x_{k(2)i,k(2)j(2)}}\\
&
=\mu_{j(1)}^{2k(1)}\delta_{k(2),-k(1)}\delta_{j(1),j(2)}=\varphi(c_{j(1)}^{k(1)}c_{j(2)}^{k(2)})
\end{align*}
\item It is easily seen that $\varphi^{\overline{\epsilon}}(\gamma_{i,j}^{k(1)}\dots
\gamma_{i,j}^{k(w)})$ is independent of $i \in \N_{*}$.
\item This is a consequence of Lemma \ref{rep}.
\item This follows from Lemma \ref{rand} almost surely.
\end{enumerate}
Thus, by Theorem \ref{TCL}, we have, almost surely, for all $p \in
\N_{*}$, $(k(1),\dots,k(p))\in \{-1,1\}^{p}$ and all
$(j(1),\dots,j(p)) \in \{1,\dots,k\}^{p}$:
$$
\lim_{n\rightarrow +\infty}\varphi^{\overline{\epsilon}}(s_{n,
j(1)}^{k(1)}\dots s_{n, j(p)}^{k(p)})=\left\{
\begin{array}{cl}
0 &\text{if $p$ is odd}\vspace{0.5 cm}\\
\displaystyle{\sum_{\substack{\mathcal{V} \in \mathcal{P}_{2}(1,\dots,2r)\\
\mathcal{V}=\{(s_{l},t_{l})_{l=1}^{l=r}\}
}}}q^{i(\mathcal{V})}\Prod_{l=1}^{r}
\varphi(c_{j(s_{l})}^{k(s_{l})}c_{j(t_{l})}^{k(t_{l})}) &
\text{if\;} p=2r
\end{array}
\right.$$ By Lemma \ref{momdef} we see that all $*$-moments of the
family $(s_{n,j})_{j \in \{1,\dots,k\}}$  converge when $n$ goes to
infinity to the corresponding  $*$-moments of the family
$(c_{j})_{j\in \{1,\dots,k\}}$:

\begin{prop} For all $p \in \N_{*}$, $(j(1),\dots,j(p)) \in
\{1,\dots,k\}^{p}$ and for all \\$(k(1),\dots,k(p))\in \{-1,1\}^{p}$
we have:
\begin{equation}\label{eq9}
\lim_{n\rightarrow +\infty}\varphi^{\overline{\epsilon}}(s_{n,
j(1)}^{k(1)}\dots s_{n, j(p)}^{k(p)})=\varphi(c_{j(1)}^{k(1)}\dots
c_{j(p)}^{k(p)})\quad\text{almost surely}
\end{equation}
\end{prop}

{\noindent \it Remark.} It is possible (and maybe easier) to apply
directly Speicher's Theorem to the independent family
$(\beta_{i,j})_{(i,j)\in \overline{I}^2}$. Then, it suffices to
follow the analogies between the Baby Fock and the $q$-Fock
frameworks to deduce the previous Proposition.

\subsection{$\Gamma_{q}(H_{\R},(U_{t})_{t \in \R})$ is QWEP}

For all $j \in \{1,\dots,k\}$ let us denote by $
g_{n,j}=\Re(s_{n,j})\quad\text{and}\quad g_{n,-j}=\Im(s_{n,j})$. By
(\ref{eq9}) we have that for all monomials $P$ in $2k$ noncommuting
variables:
\begin{equation}\label{eq6}
 \lim_{n\rightarrow
+\infty}\varphi^{\overline{\epsilon}}(P(g_{n, -k},\dots ,g_{n,
k}))=\varphi(P(G(f_{-k}),\dots, G(f_{k})))\quad\text{almost surely}
\end{equation}
Since the set of all non-commutative monomials is countable, we can
find a choice of signs $\epsilon$ such that (\ref{eq6}) is true for
all $P$. In the sequel we fix such an $\epsilon$ and forget about
the dependance on $\epsilon$.
\begin{lem}\label{mom2} For all polynomials $P$ in $2k$ noncommuting variables
we have:
\begin{equation}\label{eq7}
 \lim_{n\rightarrow
+\infty}\varphi(P(g_{n, -k},\dots ,g_{n,
k}))=\varphi(P(G(f_{-k}),\dots, G(f_{k})))
\end{equation}
\end{lem}

We are now ready to construct an embedding of
$\Gamma_{q}(H_{\R},U_{t})$ into an ultraproduct of the finite
dimensional von Neumann algebras $\Gamma_{n}$. To do so we  need to
have a uniform bound on the operators $g_{n,j}$. Let $C>0$ such that
for all $j\in I$, $\|G(f_{j})\|<C$, as in the tracial case, we
replace the $g_{n,j}$ by the their truncations
$\tilde{g}_{n,j}=\chi_{]-C,C[}(g_{n,j})g_{n,j}$. The following is
the analogue of Lemma \ref{pres}:

\begin{lem}\label{pres2} For all polynomials $P$ in $2k$ noncommuting variables
we have:
\begin{equation}\label{eq10}
 \lim_{n\rightarrow
+\infty}\varphi(P(\tilde{g}_{n, -k},\dots ,\tilde{g}_{n,
k}))=\varphi(P(G(f_{-k}),\dots, G(f_{k})))
\end{equation}
\end{lem}

{\noindent \it Remark.} For all $n \in \N_{*}$ and all $j\in I$ the
element $g_{n,j}$ is entire for the modular group (this is always
the case in a finite dimensional framework). By (3) of proposition
\ref{modu}, we have for all $j\in \{1,\dots,k\}$
$$\sigma_{z}(s_{n,j})=\lambda_{j}^{iz}s_{n,j}\quad\text{for all\;}z \in \C $$
Thus for all $z \in \C$,
\begin{equation}\label{gm}
\sigma_{z}(g_{n,j})= \left\{
\begin{array}{ll}
 \cos(z\ln(\lambda_{j}))g_{n,j}-\sin(z\ln(\lambda_{j}))g_{n,-j}&\text{for all\;}j \in\{1,\dots,k\}\\
 \sin(z\ln(\lambda_{-j}))g_{n,-j}+\cos(z\ln(\lambda_{-j}))g_{n,j}&\text{for all\;}j \in\{-1,\dots,-k\}\\
\end{array}
\right. \end{equation}

\medskip

\noindent{\it Proof of Lemma \ref{pres2}.} It suffices to show that
for all $(j(1),\dots,j(p)) \in I^{p}$ we have
$$\lim_{n\rightarrow
+\infty}\varphi(\tilde{g}_{n,j(1)}\dots \tilde{g}_{n,
j(p)})=\varphi(G(f_{j(1)})\dots G(f_{j(p)}))$$ By (\ref{eq7}) it is
sufficient to prove that
$$\lim_{n \rightarrow +\infty}|\varphi(g_{n,j(1)}\dots g_{n,
j(p)})-\varphi(\tilde{g}_{n,j(1)}\dots \tilde{g}_{n, j(p)})|=0$$
Using multi-linearity we can write
\begin{align*}
|\varphi(g_{n,j(1)}\dots g_{n,
j(p)})-\varphi(\tilde{g}_{n,j(1)}&\dots \tilde{g}_{n, j(p)})|\\
& =|\sum_{l=1}^{p}\varphi[\tilde{g}_{n,j(1)}\dots
\tilde{g}_{n,j(l-1)}(g_{n,j(l)}-\tilde{g}_{n,j(l)})g_{n,j(l+1)}\dots
g_{n,j(p)}]|\\
& \le \sum_{l=1}^{p}|\varphi[\tilde{g}_{n,j(1)}\dots
\tilde{g}_{n,j(l-1)}(g_{n,j(l)}-\tilde{g}_{n,j(l)})g_{n,j(l+1)}\dots
g_{n,j(p)}]|\\
\end{align*}
Fix $l \in \{1,\dots,p\}$, using the modular group we have:
\begin{align*}
|\varphi[\tilde{g}_{n,j(1)}\dots
\tilde{g}_{n,j(l-1)}(g_{n,j(l)}-&\tilde{g}_{n,j(l)})g_{n,j(l+1)}\dots
g_{n,j(p)}]|\\
&=|\varphi[\sigma_{i}(g_{n,j(l+1)}\dots
g_{n,j(p)})\tilde{g}_{n,j(1)}\dots
\tilde{g}_{n,j(l-1)}(g_{n,j(l)}-\tilde{g}_{n,j(l)})]|
\end{align*}
Estimating by Cauchy-Schwarz's inequality we obtain:
\begin{align*}
|\varphi[\sigma_{i}(g_{n,j(l+1)}&\dots
g_{n,j(p)})\tilde{g}_{n,j(1)}\dots
\tilde{g}_{n,j(l-1)}(g_{n,j(l)}-\tilde{g}_{n,j(l)})]|\\
&\le \varphi[\sigma_{i}(g_{n,j(l+1)}\dots
g_{n,j(p)})\tilde{g}_{n,j(1)}\dots \tilde{g}_{n,j(l-1)}^{2}\dots
\tilde{g}_{n,j(1)}\sigma_{-i}(g_{n,j(p)}\dots
g_{n,j(l+1)})]^{\f{1}{2}}\\
&\qquad \times \varphi[(g_{n,j(l)}-\tilde{g}_{n,j(l)})^{2}]^{\f{1}{2}}\\
&\le C^{l-1}\varphi[\sigma_{i}(g_{n,j(l+1)}\dots
g_{n,j(p)})\sigma_{-i}(g_{n,j(p)}\dots
g_{n,j(l+1)})]^{\f{1}{2}}\varphi[(g_{n,j(l)}-\tilde{g}_{n,j(l)})^{2}]^{\f{1}{2}}
\end{align*}
The conclusion follows from the convergence of this last term to
$0$. Indeed, by (\ref{gm}) there exists a polynomial in $2k$
non-commutative variables $Q$, independent on $n$, such that
$Q(g_{n,-k}\dots g_{n,k})=\sigma_{i}(g_{n,j(l+1)}\dots
g_{n,j(p)})\sigma_{-i}(g_{n,j(p)}\dots g_{n,j(l+1)})$. It follows by
(\ref{eq7}) that
$$\lim_{n\to +\infty}\varphi[\sigma_{i}(g_{n,j(l+1)}\dots
g_{n,j(p)})\sigma_{-i}(g_{n,j(p)}\dots
g_{n,j(l+1)})]=\varphi(Q(G(f_{-k})\dots G(f_{k}))).$$ And by Lemma
\ref{conc}, $\varphi[(g_{n,j(l)}-\tilde{g}_{n,j(l)})^{2}]$ converges
to $0$ when $n$ goes to infinity. \Eproof

\bigskip

Let us denote by $\mathcal{P}$ the $w^{*}$-dense $*$-subalgebra of
$\Gamma_{q}(H_{\R},U_{t})$ generated by the set $\{G(f_{j}),\;j\in
I\}$. We know that $\mathcal{P}$ is isomorphic to the algebra of
non-commutative polynomials in $2k$ variables (see the remark after
Lemma \ref{conc}). Given $\mathcal{U}$ a non trivial ultrafilter on
$\N$, it is thus possible to define the following $*$-homomorphism
$\Phi$ from $\mathcal{P}$ into the von Neumann ultraproduct
$\mathcal{A}=\Prod\limits_{n, \mathcal{U}}\Gamma_{n}$ by:
$$\Phi(P(G(f_{-k}),\dots,G(f_{k})))=(P(\tilde{g}_{n, -k},\dots ,\tilde{g}_{n,
k}))_{n \in \N}^{\bullet}$$  Indeed the right term is well defined
since it is uniformly bounded in norm. Let us check the hypothesis
of Theorem \ref{qwep}.
\begin{enumerate}
\item By Lemma \ref{pres2}, $\Phi$ is state
preserving.
\item It is sufficient to check that condition 2. of Lemma \ref{p} is satisfied for
every generator $\Phi (G(f_{j}))$, $j\in I$. Let us fix $j\in I$ and
recall that by (\ref{gm}) there are complex numbers $\nu_{j}$ and
$\omega_{j}$ (independent of $n$) such that
$\sigma_{-i}^{n}(g_{n,j})=\nu_{j}g_{n,j}+\omega_{j}g_{n,-j}$. We
show that condition 2. of Lemma \ref{p} is satisfied for
$x=\Phi(G(f_{j}))$ and
$z=\nu_{j}\Phi(G(f_{j}))+\omega_{j}\Phi(G(f_{-j}))$. By
$w^{*}$-density it is sufficient to consider $y=(y_{n})_{n \in
\N}^{\bullet} \in \widetilde{\mathcal{A}}$. Using Lemma \ref{pres2}
we have:
\begin{align*}
\varphi(\Phi(G(f_{j}))y)&=\lim_{n,\mathcal{U}}\varphi_{n}(\tilde{g}_{n,j}y_{n})=
\lim_{n,\mathcal{U}}\varphi_{n}(g_{n,j}y_{n})=\lim_{n,\mathcal{U}}\varphi_{n}(y_{n}\sigma_{-i}^{n}(g_{n,j}))\\
&=\lim_{n,\mathcal{U}}\varphi_{n}(y_{n}(\nu_{j}g_{n,j}+\omega_{j}g_{n,-j}))=
\lim_{n,\mathcal{U}}\varphi_{n}(y_{n}(\nu_{j}\tilde{g}_{n,j}+\omega_{j}\tilde{g}_{n,-j}))\\
&=\varphi(y(\nu_{j}\Phi(G(f_{j}))+\omega_{j}\Phi(G(f_{-j}))))
\end{align*}
\item It suffices to check that the intertwining condition given in the remark of Theorem
\ref{qwep} is satisfied for the generators $\Phi(G(f_{j}))
=(\tilde{g}_{n,j})_{n\in N}^{\bullet}$:
$$\text{for all\;} j \in I,\quad\sigma_{t}(p\Phi(G(f_{j}))p)=p\Phi(\sigma_{t}(G(f_{j})))p$$
To fix ideas we will suppose that $j\ge 0$. Recall that in this case
for all $t\in \R$ and for all $n\in \N$, we have
$$\sigma_{t}^{n}(g_{n,j})=\cos(t\ln(\lambda_{j}))g_{n,j}-\sin(t\ln(\lambda_{j}))g_{n,-j}.$$
Since the functional calculus commutes with automorphisms, for all
$t\in \R$ and for all $n\in \N$, we have:
$$\sigma_{t}^{n}(\tilde{g}_{n,j})=h(\sigma_{t}^{n}(g_{n,j})),$$
where $h(\lambda)=\chi_{]-C,C[}(\lambda)\lambda$, for all $\lambda
\in \R$. But by Lemma \ref{mom2},
$$\sigma_{t}^{n}(g_{n,j})=\cos(t\ln(\lambda_{j}))g_{n,j}-\sin(t\ln(\lambda_{j}))g_{n,-j}$$
converges in distribution  to
$$\cos(t\ln(\lambda_{j}))G(f_{j})-\sin(t\ln(\lambda_{j}))G(f_{-j})=\sigma_{t}(G(f_{j}))$$
and $\|\sigma_{t}(G(f_{j}))\|=\|G(f_{j})\|<C$. Thus, by
 Lemma \ref{conc}, we deduce that
$\sigma_{t}^{n}(\tilde{g}_{n,j})$ converges in distribution to
$\sigma_{t}(G(f_{j}))$. On the other hand, by Lemma \ref{pres2},
$$\cos(t\ln(\lambda_{j}))\tilde{g}_{n,j}-\sin(t\ln(\lambda_{j}))\tilde{g}_{n,-j}$$
also converges in distribution to
$$\cos(t\ln(\lambda_{j}))G(f_{j})-\sin(t\ln(\lambda_{j}))G(f_{-j})=\sigma_{t}(G(f_{j})).$$
Let $y \in \mathcal{A}$, using Raynaud's results we compute:
\begin{align*}
\varphi(\sigma_{t}(p\Phi(G(f_{j}))p)pyp)&=\varphi((\Delta_{n}^{it})^{\bullet}p\Phi(G(f_{j}))
p(\Delta_{n}^{-it})^{\bullet}pyp)\\
&=\varphi(p(\Delta_{n}^{it})^{\bullet}\Phi(G(f_{j}))
(\Delta_{n}^{-it})^{\bullet}pyp)\\
&=\varphi((\Delta_{n}^{it})^{\bullet}\Phi(G(f_{j}))
(\Delta_{n}^{-it})^{\bullet}py)
\end{align*}
Let $z=(z_n)_{n\in \N}^{\bullet} \in \widetilde{\mathcal{A}}.$ By
our previous observations, we have:
\begin{align*}
\varphi((\Delta_{n}^{it})^{\bullet}\Phi(G(f_{j}))
(\Delta_{n}^{-it})^{\bullet}z)&=\lim_{n,\mathcal{U}}\varphi_{n}(\Delta_{n}^{it}\tilde{g}_{n,j}\Delta_{n}^{-it}z_n)\\
&=\lim_{n,\mathcal{U}}\varphi_{n}(\sigma_{t}^{n}(\tilde{g}_{n,j})z_n)\\
&=\varphi(\sigma_{t}(G(f_{j}))z)\\
&=\lim_{n,\mathcal{U}}
\varphi_{n}((\cos(t\ln(\lambda_{j}))\tilde{g}_{n,j}-\sin(t\ln(\lambda_{j}))\tilde{g}_{n,-j})z_n)\\
&=\varphi((\cos(t\ln(\lambda_{j}))\Phi(G(f_{j}))-\sin(t\ln(\lambda_{j}))\Phi(G(f_{-j})))z)\\
&=\varphi((p\Phi(\sigma_{t}(G(f_{j})))p)zp)\\
\end{align*}
By $w^{*}$-density and continuity, we can replace $z$ by $py$ in the
previous equality, which gives:
$$\varphi(\sigma_{t}(p\Phi(G(f_{j}))p)pyp)=\varphi((p\Phi(\sigma_{t}(G(f_{j})))p)pyp). $$
Thus, taking
$y=\sigma_{t}(p\Phi(G(f_{j}))p)-p\Phi(\sigma_{t}(G(f_{j})))p$, and
by the faithfulness of $\varphi(p\,.\,p)$ we deduce that
$$\sigma_{t}(p\Phi(G(f_{j}))p)=p\Phi(\sigma_{t}(G(f_{j})))p \in
p\text{Im}(\Phi)p$$
\end{enumerate}

By Theorem \ref{qwep}, $\Theta=p\Phi p$ can be extended into a
(necessarily injective because state preserving) $w^{*}$-continuous
$*$-homomorphism from $\Gamma_{q}(H_{\R},U_{t})$ into
$p\mathcal{A}p$ with a completely complemented image. By its
corollary \ref{corqwep}, since the algebras $\Gamma_{n}$ are finite
dimensional and a fortiori are QWEP, it follows that
$\Gamma_{q}(H_{\R},U_{t})$ is QWEP.

\begin{thm}\label{qwepf}
If $H_{\R}$ is a finite dimensional real Hilbert space equipped with
a group of orthogonal transformations $(U_{t})_{t \in \R}$, then the
von Neumann algebra $\Gamma_{q}(H_{\R},U_{t})$ is QWEP.
\end{thm}

\noindent{\it Remark.} We have only proved the Theorem for $H_{\R}$
of even dimension over $\R$. We did this only for simplicity of
notations. Of course this is not relevant since, if the dimension of
$H_{\R}$ is odd, then we just have to consider the real Hilbert
space $H_{\R}\oplus \R$ equipped with $(U_{t}\oplus\text{Id})_{t \in
\R}$. $\Gamma_{q}(H_{\R}\oplus \R,U_{t}\oplus\text{Id})$ is QWEP by
our previous discussion. Let us denote by $Q$ the projection from
$H_{\R}\oplus \R$ onto $H_{\R}$, then $Q$ intertwines
$(U_{t}\oplus\text{Id})_{t \in \R}$ and $(U_{t})_{t \in \R}$. In
this situation we can consider $\Gamma_{q}(Q)$, the second
quantization of $Q$ (see \cite{Hi}), which is a conditional
expectation from $\Gamma_{q}(H_{\R}\oplus \R,U_{t}\oplus\text{Id})$
onto $\Gamma_{q}(H_{\R},U_{t})$. Thus $\Gamma_{q}(H_{\R},U_{t})$ is
completely complemented into a QWEP von Neumann algebra, so
$\Gamma_{q}(H_{\R},U_{t})$ is QWEP.

\begin{cor}\label{qwepa} If $(U_{t})_{t \in \R}$  is almost periodic
on $H_{\R}$, then $\Gamma_{q}(H_{\R},U_{t})$ is QWEP.
\end{cor}

\proof There exist an invariant real Hilbert space $H_{1}$, an
orthogonal family of invariant $2$-dimensional real Hilbert spaces
$(H_{\alpha})_{\alpha \in A}$ and real eigenvalues
$(\lambda_{\alpha})_{\alpha \in A}$ greater than $1$ such that
$$H_{\R}=H_{1} \mathop{\oplus}_{\alpha \in A}H_{\alpha}\quad\text{and}
\quad U_{t|H_{1}}=\text{Id}_{H_{1}},\quad U_{t|H_{\alpha}}= \left(
\begin{array}{cc}
\cos(t\ln(\lambda_{\alpha})) & -\sin(t\ln(\lambda_{\alpha}))\\
\sin(t\ln(\lambda_{\alpha})) & \cos(t\ln(\lambda_{\alpha}))
\end{array}
\right)$$ In particular it is possible to find a net
$(I_{\beta})_{\beta \in B}$ of isometries from finite dimensional
subspaces $H_{\beta}\subset H_{\R}$ into $H_{\R}$, such that for all
$\beta \in B$, $H_{\beta}$ is stable by $(U_{t})_{t \in \R}$ and
$\bigcup\limits_{\beta \in B}H_{\beta}$ is dense in $H_{\R}$. By
second quantization, for all $\beta \in B$, there exists an
isometric $*$-homomorphism $\Gamma_{q}(I_{\beta})$ from
$\Gamma_{q}(H_{\beta},U_{t|H_{\beta}})$ into
$\Gamma_{q}(H_{\R},U_{t})$, and $\Gamma_{q}(H_{\R},U_{t})$ is the
inductive limit (in the von Neumann algebra's sense) of the algebras
$\Gamma_{q}(H_{\beta},U_{t|H_{\beta}})$. By the previous Theorem,
for all $\beta \in B$, $\Gamma_{q}(H_{\beta},U_{t|H_{\beta}})$ is
QWEP, thus $\Gamma_{q}(H_{\R},U_{t})$ is QWEP, as an inductive limit
 of QWEP von Neumann algebras. \Eproof

\section{The general case}\label{section5}

We will derive the general case by discretization and an
ultraproduct argument similar to that of the previous section.

\subsection{Discretization argument}

Let $H_{\R}$ be a real Hilbert space and $(U_{t})_{t \in \R}$ a
strongly continuous group of orthogonal transformations on $H_{\R}$.
We denote by $H_{\C}$ the complexification of $H_{\R}$ and by
$(U_{t})_{t \in \R}$ its extension to a group of unitaries on
$H_{\C}$. Let $A$ be the (unbounded) non degenerate positive
infinitesimal generator of $(U_{t})_{t \in \R}$. For every $n \in
\N_{*}$ let $g_n$ be the bounded Borelian function defined by:
$$g_n=\chi_{]1,1+\f{1}{2^n}[}+\left(\sum_{k=2^{n}+1}^{n2^{n}-1}\f{k}{2^n}
\chi_{[\f{k}{2^n},\f{k+1}{2^n}[}\right)+n\chi_{[n,+\infty[}$$ and
$$f_{n}(t)=g_{n}(t)\chi_{\{t>1\}}(t)+\f{1}{g_{n}(1/t)}\chi_{\{t<1\}}(t)+\chi_{\{1\}}(t)\quad\text{for
all\;} t\in \R_{+}$$ It is clear that
\begin{equation}\label{inv}
f_{n}(t)\nearrow t\quad\text{for all}\quad t\ge
1\qquad\text{and}\qquad f_{n}(t)=\f{1}{f_{n}(1/t)}\quad\text{for
all}\quad t\in \R_{+}^{*}.
\end{equation}
 For all $n \in
\N_{*}$, let $A_n$ be the invertible positive and bounded operator
on $H_{\C}$ defined by $A_n=f_{n}(A)$. Denoting by $\mathcal{J}$ the
conjugation on $H_{\C}$, we know, by \cite{Sh}, that
$\mathcal{J}A=A^{-1}\mathcal{J}$. By the second part of (\ref{inv}),
it follows that for all $n\in \N_{*}$,
\begin{equation}\label{inva}
\mathcal{J}A_n =\mathcal{J}f_{n}(A)=f_{n}(A^{-1})\mathcal{J}=
f_{n}(A)^{-1}\mathcal{J}=A_{n}^{-1}\mathcal{J}
\end{equation}
Consider the strongly continuous unitary group $(U_{t}^{n})_{t\in
\R}$ on $H_{\C}$ with positive non degenerate and bounded
infinitesimal generator given by $A_{n}$. By definition, we have
$U_{t}^{n}=A_{n}^{it}$. By (\ref{inva}), and since $\mathcal{J}$ is
anti-linear, we have for all $n \in \N_{*}$ and all $t \in \R$:
$$\mathcal{J}U_{t}^{n}=\mathcal{J}A_{n}^{it}=A_{n}^{it}\mathcal{J}
=U_{t}^{n}\mathcal{J}$$ It follows that for all $n \in \N_{*}$ and
for all $t \in \R$, $H_{\R}$ is globally invariant by $U_{t}^{n}$,
thus we have
$$U_{t}^{n}(H_{\R})=H_{\R}$$
Hence,  $(U_{t}^{n})_{t\in \R}$ induces a group of orthogonal
transformations on $H_{\R}$ such that its extension on $H_{\C}$ has
infinitesimal generator given by the discretized operator  $A_{n}$.
In the following we will index by $n \in \N_{*}$ the objects
relative to the discretized von Neumann algebra
$\Gamma_{n}=\Gamma_{q}\left(H_{\R},(U_{t}^{n})_{t\in \R}\right)$. We
simply set $\Gamma=\Gamma_{q}\left(H_{\R},(U_{t})_{t\in\R}\right)$.

\medskip

\noindent{\it Remark.} Notice that $H_{\C}$ is contractively
included in $H$ and all $H_{n}$, and that the inclusion
$H_{\R}\subset H$ (respectively $H_{\R}\subset H_n$) is isometric
since $\Re(\scal{\,.\,}{\,.\,}_{U})_{\left|H_{\R}\times
H_{\R}\right.}=\scal{\,.\,}{\,.\,}_{H_{\R}}$ (see \cite{Sh}).
Moreover for all $n \in \N_{*}$ the scalar products
$\scal{\,.\,}{\,.\,}_{U^{n}}$ and $\scal{\,.\,}{\,.\,}_{H_{\C}}$ are
equivalent on $H_{\C}$ since $A_{n}$ is bounded.
\begin{sch}\label{con} For all $\xi$ and $\eta$ in $H_{\C}$ we have:
$$\lim_{n\to +\infty}\scal{\xi}{\eta}_{H_{n}}=\scal{\xi}{\eta}_{H}$$
\end{sch}

\proof Let $E_{A}$ be the spectral resolution of $A$. Take $\xi \in
H_{\C}$ and denote by $\mu_{\xi}$ the finite positive measure on
$\R_{+}$ given by $\mu_{\xi}=\scal{E_{A}(.)\xi}{\xi}_{H_{\C}}$.
Since for all $\lambda \in \R_{+}$, $\lim\limits_{n \to
+\infty}g\circ f_{n}(\lambda)=g(\lambda)$, and
$g(\lambda)=2\lambda/(1+\lambda)$ is bounded on $\R_{+}$, we have by
the Lebesgue dominated convergence Theorem:
\begin{align*}
\|\xi\|^{2}_{H}&=\scal{\f{2A}{1+A}\xi}{\xi}_{H_{\C}}=\int_{\R_{+}}g(\lambda)
\dd \mu_{\xi}(\lambda)\\
&=\lim_{n \to +\infty}\int_{\R_{+}}g\circ f_{n}(\lambda) \dd
\mu_{\xi}(\lambda)=\lim_{n \to
+\infty}\scal{\f{2A_{n}}{1+A_{n}}\xi}{\xi}_{H_{\C}}=\lim_{n \to
+\infty}\|\xi\|^{2}_{H_{n}}
\end{align*}
And we finish the proof by polarization.\Eproof

Let $E$ be the vector space given by
$$E=\cup_{k\in \N_{*}}\chi_{[\f{1}{k},k]}(A)(H_{\R})$$
We have
$$\mathcal{J}\chi_{[\f{1}{k},k]}(A)=\chi_{[\f{1}{k},k]}(A^{-1})\mathcal{J}
=\chi_{[\f{1}{k},k]}(A)\mathcal{J}$$ thus $E\subset H_{\R}$.  Since
$A$ is non degenerate,
$$\overline{\cup_{k\in
\N_{*}}\chi_{[\f{1}{k},k]}(A)(H_{\C})}=\chi_{]0,+\infty[}(A)(H_{\C})=H_{\C}$$
It follows that $E$ is dense in $H_{\R}$. Let $(e_{i})_{i \in I}$ be
an algebraic basis of unit vectors of $E$ and denote by
$\mathcal{E}$ the algebra generated by the Gaussians $G(e_{i})$ for
$i \in I$. $\mathcal{E}$ is $w^{*}$ dense in $\Gamma$ and every
element in $\mathcal{E}$ is entire for $(\sigma_{t})_{t \in \R}$
(because for all $k \in \N_{*}$, $A$ is bounded and has a bounded
inverse on $\chi_{[\f{1}{k},k]}(A)(H_{\C})$). Denoting by $W$ the
Wick product in $\Gamma$, we have for all $i \in I$ and all $z \in
\C$:
\begin{equation}\label{ent}
\sigma_{z}(G(e_{i}))=W(U_{-z}e_{i})=W(A^{-iz}e_{i})
\end{equation}
Since $H_{\R}\subset H$ and for all $n\in \N_{*}$, $H_{\R}\subset
H_n$(isometrically), we have by (\ref{moman})
\begin{equation}\label{con3}\text{For all}\quad(i,n)\in
I\times\N_{*},\qquad\|G_{n}(e_{i})\|=\f{2}{\sqrt{1-q}}
\end{equation}
\begin{sch}\label{con2} For all $r\in \R$ and for all $i \in I$ we
have
$$\sup_{n\in \N_{*}}\|\sigma_{ir}^{n}(G_{n}(e_{i}))\|<+\infty$$
\end{sch}

\proof Fix $i \in I$. By (\ref{ent}):
\begin{align*}
\|\sigma_{ir}^{n}(G_{n}(e_{i}))\|&=\|W(A_{n}^{r}e_{i})\|=
\|a_{n}^{*}(A_{n}^{r}e_{i})+a_{n}(\mathcal{J}A_{n}^{r}e_{i})\|\\
&\le
C_{|q|}^{\f{1}{2}}(\|A_{n}^{r}e_{i}\|_{H_{n}}+\|\mathcal{J}A_{n}^{r}e_{i}\|_{H_{n}})\\
&\le
C_{|q|}^{\f{1}{2}}(\|A_{n}^{r}e_{i}\|_{H_{n}}+\|\Delta_{n}^{\f{1}{2}}A_{n}^{r}e_{i}\|_{H_{n}})\\
&\le
C_{|q|}^{\f{1}{2}}(\|A_{n}^{r}e_{i}\|_{H_{n}}+\|A_{n}^{r-\f{1}{2}}e_{i}\|_{H_{n}})
\end{align*}
Thus it suffices to prove that for all $r\in \R$ we have
$$\sup_{n\in \N_{*}}\|A_{n}^{r}e_{i}\|_{H_{n}}<+\infty$$
Let us denote by $\mu_{i}=\scal{E_{A}(.)e_{i}}{e_{i}}_{H_{\C}}$ and
by $g_{r}(\lambda)=2\lambda^{2r+1}/(1+\lambda)$. There exists $k \in
\N_{*}$ such that $e_{i}\in \chi_{[1/k,k]}(A)(H_{\R})$, thus we have
:$$\|A_{n}^{r}e_{i}\|_{H_{n}}^{2}=\scal{g_{r}\circ
f_{n}(A)e_{i}}{e_{i}}_{H_{\C}}=\int_{[1/k,k]}g_{r}\circ
f_{n}(\lambda)\dd \mu_{i}(\lambda)$$ It is easily seen that
$(g_{r}\circ f_{n})_{n \in \N_{*}}$ converges uniformly to $g_{r}$
on $[1/k,k]$. The result follows by:
$$\lim_{n \to +\infty}\|A_{n}^{r}e_{i}\|_{H_{n}}^{2}=\lim_{n \to +\infty}\int_{[1/k,k]}g_{r}\circ
f_{n}(\lambda)\dd \mu_{i}(\lambda)=\int_{[1/k,k]}g_{r}(\lambda)\dd
\mu_{i}(\lambda)=\|A^{r}e_{i}\|_{H}^{2}.$$ \Eproof

\subsection{Conclusion}

Recall that $\mathcal{E}$ is isomorphic to the complex free
$*$-algebra with $|I|$ generators. Let $\mathcal{U}$ be a free
ultrafilter on $\N_{*}$, by (\ref{con3}) we can define a
$*$-homomorphism $\Phi$ from $\mathcal{E}$ into the von Neumann
algebra ultraproduct over $\mathcal{U}$ of the algebras $\Gamma_{n}$
by:
\begin{equation*}
\begin{split}
\Phi\;:\;\mathcal{E}&\longrightarrow \mathcal{A}=\Prod_{n, \mathcal{U}}\Gamma_{n}\\
G(e_{i})&\longmapsto (G_{n}(e_{i}))_{n\in \N_{*}}^{\bullet}
\end{split}
\end{equation*}
We will now check the hypothesis of Theorem \ref{qwep}.
\begin{enumerate}
\item We first check that $\Phi$ is state preserving. It suffices
to verify it for a product of an even number of Gaussians. Take
$(i_{1},\dots,i_{2k})\in I^{2k}$, we have by Scholie \ref{con}:
\begin{align*}
\varphi(G(e_{i_{1}})\dots
G(e_{i_{2k}}))&=\sum_{\substack{\mathcal{V}\in
\mathcal{P}_{2}(1,\dots,k)\\
\mathcal{V}=((s(l),t(l)))_{l=1}^{l=k}}}
q^{i(\mathcal{V})}\Prod_{l=1}^{l=k}
\scal{e_{i_{s(l)}}}{e_{i_{t(l)}}}_{H}\\
&=\lim_{n \to +\infty}\sum_{\substack{\mathcal{V}\in
\mathcal{P}_{2}(1,\dots,k)\\
\mathcal{V}=((s(l),t(l)))_{l=1}^{l=k}}}
q^{i(\mathcal{V})}\Prod_{l=1}^{l=k}
\scal{e_{i_{s(l)}}}{e_{i_{t(l)}}}_{H_{n}}\\
&=\lim_{n \to +\infty}\varphi_{n}(G_{n}(e_{i_{1}})\dots
G_{n}(e_{i_{2k}}))
\end{align*}
This implies, in particular that $\Phi$ is state preserving.
\item Condition 1. of lemma \ref{p} is satisfied by Scholie \ref{con2}.
\item It suffices to check that for all $i\in I$ and all $t \in \R$,
 $(\sigma_{t}^{n}(G_{n}(e_{i})))_{n\in \N_{*}}^{\bullet}
\in \overline{\text{Im}\Phi}^{w^{*}}.$ Fix $i \in I$ and $t \in \R$.
For all $n \in \N_{*}$ we have
$$\|A_{n}^{-it}e_{i}-A^{-it}e_{i}\|_{H_{\R}}^{2}=\int_{\R_{+}}|f_{n}^{-it}(\lambda)-\lambda^{-it}|^{2}\dd
\mu_{i}(\lambda)$$ By the Lebesgue dominated convergence Theorem, it
follows that
$$\lim_{n\to +\infty}\|A_{n}^{-it}e_{i}-A^{-it}e_{i}\|_{H_{\R}}=0.$$
By (\ref{con3}) we deduce that
$$\lim_{n\to
+\infty}\|G_{n}(A_{n}^{-it}e_{i})-G_{n}(A^{-it}e_{i})\|=0$$ Thus we
have
$$(\sigma_{t}^{n}(G_{n}(e_{i})))_{n\in \N_{*}}^{\bullet}=(G_{n}(A_{n}^{-it}e_{i}))_{n\in \N_{*}}^{\bullet}
=(G_{n}(A^{-it}e_{i}))_{n\in \N_{*}}^{\bullet}\in
\overline{\text{Im}\Phi}^{\|.\|}\subset
\overline{\text{Im}\Phi}^{w^{*}}.$$
\end{enumerate}

By Theorem \ref{qwep}, we deduce our main Theorem:

\begin{thm}\label{qwepc} Let $H_{\R}$ be a real Hilbert space given with a group
of orthogonal transformations $(U_{t})_{t \in \R}$. Then for all $q
\in (-1,1)$ the $q$-Araki-Woods algebra
$\Gamma_{q}(H_{\R},(U_{t})_{t \in \R})$ is QWEP.
\end{thm}

\bigskip

{\noindent \it Remark.} We were unable to prove that the
$C^{*}$-algebra $C^{*}_{q}(H_{\R},(U_{t})_{t \in \R})$ (for
$(U_{t})_{t\in \R}$ non trivial) is QWEP, even for a finite
dimensional Hilbert space $H_{\R}$. The proof of Lemma \ref{wcp}
could not be directly adapted to this case. Indeed, in the
non-tracial framework, the ultracontractivity of the
$q$-Ornstein-Uhlenbeck semi-group $(\Phi_{t})_{t\in \R_{+}}$ is
known when $A$ is bounded and $t>\frac{\ln(\|A\|)}{2}$ but in any
cases it fails for $0<t<\frac{\ln(\|A\|)}{4}$ (see \cite{Hi} Theorem
4.1 and Proposition 4.5).

\addcontentsline{toc}{chapter}{Bibliography}
\bibliographystyle{plain}

\end{document}